\documentclass[1p]{elsarticle}

\journal{Journal of Mathematical Bioscience}
\usepackage{epstopdf}
\usepackage{amsmath} 
\usepackage{amssymb}  
\usepackage{tabularx}
\usepackage{amsmath}
\usepackage{amsthm}
\usepackage{float}
\usepackage{soul}
\usepackage{color}
\usepackage{rotating,booktabs,multirow}
\usepackage{graphicx}
\usepackage{hyperref}

\newtheorem{asu}{Assumption}
\newtheorem{remark}{Remark}[section]
\newtheorem{theorem}{Theorem}

\newcounter{subassumption}[asu]

\makeatletter
\renewcommand{\p@subassumption}{\theasu}
\makeatother

\begin{document}

\begin{frontmatter}

\title{Immune Therapeutic Strategies Using Optimal Controls with $L^1$ and $L^2$ Type Objectives\tnoteref{mytitnote}}

\author[eecs1]{O. Bara\corref{cor1}}
\cortext[cor1]{Corresponding author.}
\ead{obara@vols.utk.edu}

\author[eecs1]{S.M. Djouadi}
\ead{mdjouadi@utk.edu}

\author[math]{J.D Day}
\ead{judyday@utk.edu}

\author[math]{S. Lenhart}
\ead{lenhart@math.utk.edu}

\address[eecs1]{Department of Electrical Engineering  and Computer Science, University of Tennessee, Knoxville, TN 37996, United States}
\address[math]{Department of Mathematics, University of Tennessee, Knoxville, TN 37996, United States }

\tnotetext[mytitnote]{Work partially supported by the NSF-DMS Award 1122462  and a Fullbright Scholarship.}

\begin{abstract}
	Therapeutic strategies to correct an excessive immune response to pathogenic infection is investigated as an optimal control problem. The control problem is formulated around a  four dimensional mathematical model describing the inflammatory response to a pathogenic insult with two therapeutic control inputs  which have either a direct pro- or anti-inflammatory effect in the given system. We use Pontryagin's maximum principle and discuss  necessary optimality conditions. We consider both an $L^1$ type objective functional as well as an $L^2$ type objective. For the former, the presence of singular control will be addressed.  For each case, numerical simulations using a nonlinear programming optimization solver to acquire different drug treatment strategies are presented and discussed. The  results provide insight for possible treatment strategies and the methods could be a relevant tool for future practice to assist in better prediction of clinical outcomes and subsequently better treatment for patients. 
\end{abstract}


\end{frontmatter}


\section{Introduction}
The immune response is a complex mechanism triggered in an organism as a result of biological or physical stress, such as the presence of a foreign body (pathogenic infection) or trauma. The natural behavior of an organism to respond to these phenomena is through stress adaptation by trying to eliminate the invading pathogen threat in case of an infection while also promoting tissue repair due to the self-harming effects of inflammation. However, an excessive and dysregulated inflammatory response may lead to systemic inflammatory response syndrome (SIRS) and sepsis, resulting in tissue damage, organ dysfunction or even death. According to \cite{cohen2002immunopathogenesis}, the overall mortality is approximately $30\%$, rising to $40\%$ in the elderly and is $50\%$ or greater in patients with severe SIRS and sepsis.

With the withdrawal of previously approved drugs to treat sepsis (\cite{FDAXigris}) and the acknowledgment that potential therapies for sepsis such as anti-tumor necrosis factor-$\alpha$ appear to require accurate timing (\cite{Grau1997}, \cite{Reinhart2001}), identifying effective therapies and patient-specific protocols for curbing sepsis is increasingly important. Using computational tools and methodologies is one means by which to explore this and offer insight. In particular, applying automatic control for biomedical therapeutic intervention has recently gained interest in such areas as glucose control for diabetes or in critically ill patients and for anesthesia depth control, as seen for example, in
\cite{parker2010systems}, \cite{bequette2012challenges}, \cite{doyle2014closed}, \cite{haddad2003}, and the references therein.   The use of optimal control theory for biomedical applications can be seen in \cite{stengel2004stochastic}, \cite{kirschner1997}, and \cite{stengel2002optimal} for example; and, more recently, the use of model predictive control (MPC) in \cite{day2010using} and \cite{zitelli2015combining} which used the same mathematical model as used herein. MPC was also used in  \cite{parker2000advanced} and \cite{grosman2010zone} for diabetes treatment strategies and \cite{li2013} for an application to a multi-compartment respiratory system.

In the cases where MPC is used, an optimal control problem is essentially solved for each specified time interval over the duration of the time span, as though implementing the control measures in real time (i.e. online) as the response is evolving. In this work, we investigate the \emph{offline} strategy that solves the control problem for the entire time period of interest:  from the start of treatment till the end of the observation period. Within this  strategy, we explore differences in the form of the objective function used as well as determine the effectiveness of the identified treatment protocol.

As in our prior work (\cite{bara2015}, \cite{baraoptL1}, \cite{day2010using}, \cite{zitelli2015combining}), we utilize a phenomenological ordinary differential equations model of the system inflammatory response to pathogenic infection developed in \cite{reynolds2006reduced} (see also \cite{day2006reduced}). This model provides a dynamical system with rich behavior that is ideal for testing various theoretical control strategies. A summary of its important features are as follows:
\begin{itemize}
\item This model is based on the early non-specific protective mechanism, namely, the \emph{innate immune response}, in contrast to the mechanisms of \emph{adaptive immunity}. It is a lumped variable/lumped parameter model, representing an abstract view of the dynamics of an acute inflammatory response and implicitly considering affects of typical intervention strategies such as antibiotic use or the administration of resuscitating fluids;
\item An anti-inflammatory mediator is included as a dynamic variable of the model. This mediator plays an important role in directly mitigating inflammation as a way to prevent excessive tissue damage caused by severe inflammation; and
\item The model's biological relevance has been confirmed through its good qualitative reproduction of severe systemic inflammatory states as observed in a clinical setting.
\end{itemize}

In this paper, we explore the use of optimal control theory to derive an immune control strategy for the chosen model consisting  of manipulating the pro- and/or anti-inflammatory mediator  in order to reestablish the healthy equilibrium of the virtual patient.  Note that the model does not target a specific pathogen and the control doses represent concentrations of inflammatory mediators with idealized effect on the system. Even so, the careful balance required for effective treatment is evident. For instance, a large dose of pro-inflammatory therapy, while aiming to eliminate the pathogen, would also pose a risk in causing possibly irreparable damage to organs. Similarly, using a large amount of anti-inflammatory therapy could suppress the negative effects of the inflammatory-damage feedback loop and consequently preserve organ health; however, this intervention may also suppress the positive effects of inflammation in controlling pathogen which may then be uninhibited to rapidly grow. The dual objectives to eliminate pathogen but also minimize damage via the administration of pro- and anti-inflammatory mediators, respectively, necessitates the use of a control algorithm in order to find a suitable balance that produces the desired outcome: a healthy patient.\\

\indent   In \cite{baraoptL1} and \cite{bara2015} the necessary conditions of Pontryagin's maximum principle \cite{pontryagin1987} were used in order to formulate a  two point boundary value  problem (TPBVP) which was solved numerically with the forward-backward sweep method \cite{lenhart2007optimal} (indirect method). In the current work, we choose to use a direct approach to the numerical solutions that relies on approximating the original optimal control problem by first discretizing the state and control history and then solving it using nonlinear programming (NL). Further details about the solver used will be given in section \ref{sec:numeric}. Constraints are handled more easily in the direct method than in the indirect method where one must derive adjoint equations along with the transversality and optimality conditions before the TPBVP problem can be solved. Furthermore, in the indirect method, the presence of state inequality constraints requires \emph{a priori} knowledge of the number and duration of constrained arcs which is unknown beforehand \cite{bryson}.

\indent We compare the results of using a quadratic or $L^2$ objective function versus a linear or $L^1$ one. The  use of a quadratic  cost to derive the optimal control can be explained by two important facts. First the corresponding Hamiltonian will be strictly convex in the control with a unique minimum, and second, the mathematical problem is more straightforward to solve. Alternatively, using an $L^1$ or linear objective avoids distortion due to the square used in the $L^2$ form which puts a small penalty on lower doses and therefore, potentially biases the real effect of the control in the cost function \cite{ledzewicz2007optimal}. Also,  the use of an $L^1$ objective often results in ``bang-bang'' controls, meaning the control is applied at either its maximum level or its minimum level, which provide simple dosing protocols whose amounts do not vary greatly over time. In this respect, such controls are more analogous to how medical treatments are administered in practice and are considered natural choices as candidates for optimality \cite{swierniak2003optimal}.

\indent The simulation results presented herein show some similarities when compared to \cite{baraoptL1} and \cite{bara2015}, in that the optimal control doses follow  a certain strategy whether considering a  quadratic or a linear objective. As will be noted later, the optimal inputs for a virtual patient whose inflammatory mediators are elevated above some prescribed threshold are characterized by a bolus of pro-inflammatory therapy followed just after by anti-inflammatory therapy, consistent with earlier findings in \cite{day2010using}.

\indent
The paper is organized as follows. Sections \ref{sec:model} and \ref{sec:optproblem} respectively explain the mathematical model and the optimal control problem. A brief background regarding the numerical solver (PSOPT) is provided in Section \ref{sec:numeric} together with the numerical results for both the quadratic and linear cases. Suggestions for future extensions to this work are discussed in section \ref{sec:conclusion}.

\section{Immune response model} \label{sec:model}
We consider the mathematical model describing the inflammatory response for a host-pathogen interaction as proposed in \cite{reynolds2006reduced} and consider the model as representing the immune dynamics of a virtual patient. There are four system states of the model:

\begin{itemize}
\item $P$ : the bacterial pathogen population that initiates the inflammatory response;
\item $N$ : combined collection of early pro-inflammatory mediators such as activated phagocytes (e.g. neutrophils) and the pro-inflammatory cytokines they produce;
\item $D$ : a marker of global tissue damage which incites further inflammation and which is useful for determining response outcomes; and
\item $C_a$: combined collection of anti-inflammatory mediators, such as cortisol, Interleukin-10 (IL-10) and Transforming Growth Factor-beta (TGF-$\beta$).
\end{itemize}

In \cite{day2010using} it was first proposed to include non-negative, time-varying inputs to account for therapeutic controls in this model. Here, we denote these as follows:
\begin{itemize}
\item $u_p(t)$: a pro-inflammatory enhancer that provides direct input into the $N$ equation; and
\item $u_a(t)$: an anti-inflammatory enhancer that provides direct input into the $C_a$ equation.
\end{itemize}

The dynamic equations describing the interactions of the four system variables, $P$, $N$, $D$, and $C_a$, are governed by various rate constants or parameters whose numerical values can be found in \cite{reynolds2006reduced}. As was done in \cite{day2010using}, we consider different sets of parameter values to represent different virtual patients' inflammatory responses. We refer to the values given in \cite{reynolds2006reduced} as the \textit{reference set} of parameter values.

\begin{figure}[!ht]
\begin{center}
	\includegraphics[scale=0.3]{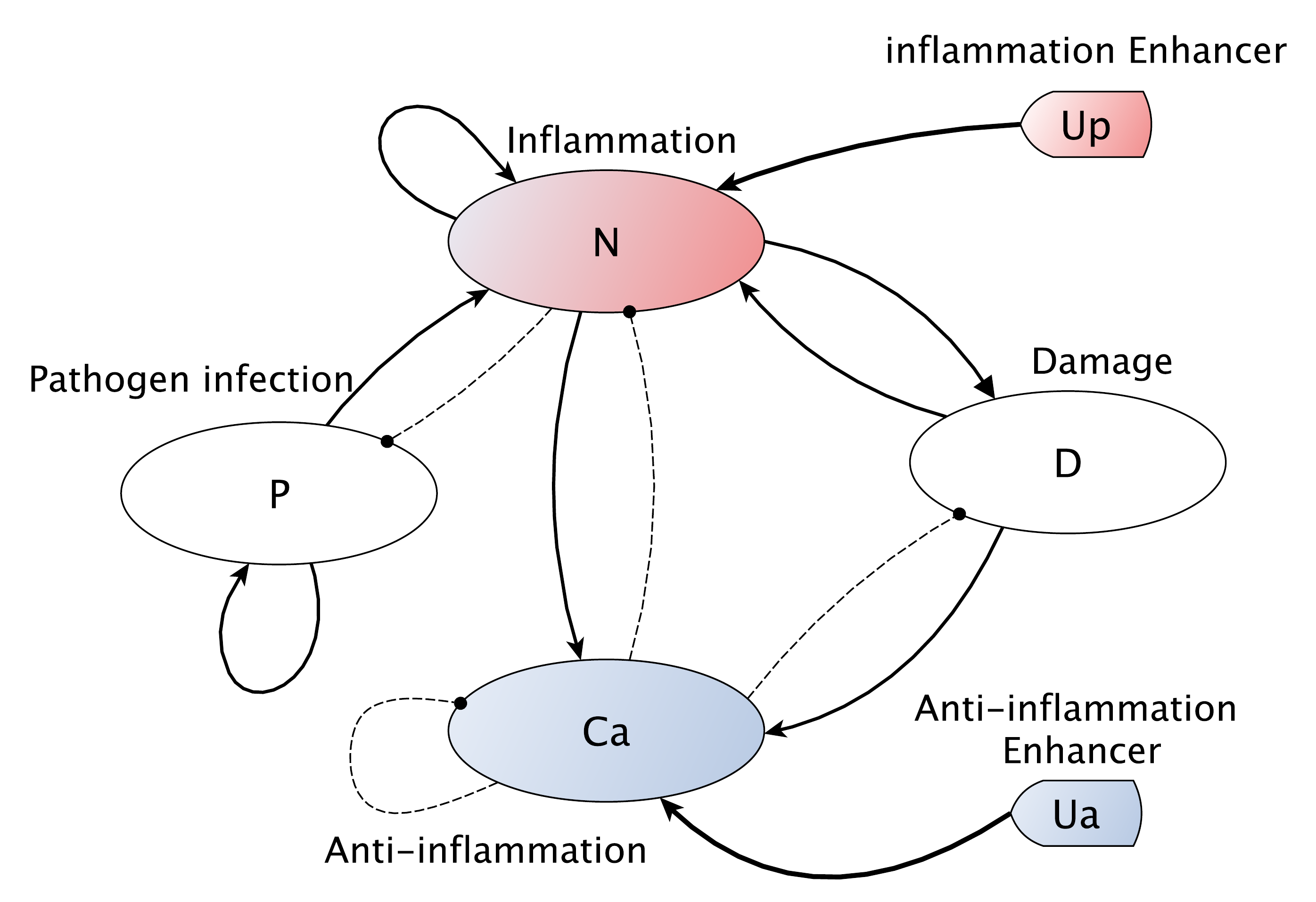}
\end{center}
\caption{Diagram of natural and enhanced immune dynamics in response to pathogenic stimulus. $P$: replicating pathogen, $N$: early pro-inflammatory immune mediators, $D$: marker of tissue damage/dysfunction caused by inflammatory response, $C_a$: inhibitory anti-inflammatory mediators, $u_a$ and $u_p$: time-varying input controls for the anti- and pro-inflammatory therapy, respectively. Solid lines with arrow heads and dashed lines with nodes/circular heads represent upregulation and inhibition, respectively.} \label{Fig:diagram}
\label{Fig:fig1}
\end{figure}

The diagram in Figure \ref{Fig:diagram} characterizes the different interactions that exist between the states of the inflammatory model. A solid line with an arrow head indicates an up-regulation, whereas a dashed line with circular head indicates inhibition or downregulation of a process. For instance, early pro-inflammatory mediators, $N$, respond to the presence of pathogen, $P$, by initiating self-recruitment of additional inflammatory mediators and $N$ is therefore upregulated by the interaction with $P$ to attempt to efficiently eliminate the pathogen. The self up-regulation that exists for $P$ is due to replication. Furthermore, $N$ suppresses $P$ by eliminating it at some rate when they interact; however, the inflammation caused by $N$ results in tissue damage, $D$, which can provide a positive feedback into the early inflammatory mediators depending on their intensity. To balance this, anti-inflammatory mediators (e.g., cortisol, IL-10, TGF-$\beta$) can mitigate the inflammation and its harmful effect by suppressing the response by $N$ and the effects of $D$ in various ways. The ordinary differential equations modeling these interactions is given by Equations (\ref{eq:mod1}) - (\ref{eq:mod4}):

\begin{align}
\dot{P}(t) &= k_{pg} P(t) \left(1-\frac{P(t)}{p_{\infty}}\right)-\frac{k_{pm} s_{m}P(t)}{\mu_{m}+k_{mp} P(t)}-k_{pn}F(N(t))P(t)\label{eq:mod1}\\
\dot{N}(t) &= \frac{s_{nr}R(P(t),N(t),D(t))}{\mu_{nr}+R(P(t),N(t),D(t))} - \mu_n N(t)+ u_p(t) \label{eq:mod2}\\
\dot{D}(t) &= k_{dn} \frac{F(N(t))^6}{x_{dn}^6 + F(N(t))^6}-\mu_{d} D(t) \label{eq:mod3}\\
\dot{C_{a}}(t) &= s_{c}+k_{cn} \frac{F(N(t)+k_{cnd}D(t))}{1+F(N(t)+k_{cnd}D(t))}-\mu_{c}C_{a}(t)+u_a(t)\label{eq:mod4},
\end{align}

\noindent where
\begin{align*}
R(P, N, D) &= F(k_{np}P(t)+k_{nn}N(t)+k_{nd}D(t)) \ \ \textrm{ and}\\[2mm]
F(x) &= \frac{x}{1+\left(\frac{C_a(t)}{c_\infty}\right)^2}.
\end{align*}

 As explained in  \cite{reynolds2006reduced}, $F(x)$ represent a Hill function of order $6$ that models the impact of activated phagocytes and their by-products ($N$) on the creation of damaged tissue. Tissue damage ($D$) increases in a switch-like sigmoidal fashion as $N$ increases, so that it takes sufficiently high levels of $N$ to incite a moderate increase in damage and that the increase in damage saturates with sufficiently elevated and sustained $N$ levels. The coefficient (exponent) 6 was chosen, therefore, to model this aspect which ensured that the healthy equilibrium had a reasonable basin of attraction for the $N/D$ system. 

The model was formulated to represent a most abstract form of the complex processes involved in the systemic inflammatory response.  Hence, the variables $N$ and $C_a$ represent multiple mediators with similar inflammatory characteristics, and D is an abstract representation of tissue damage. This abstraction reduces the description to four essential variables which also allows for tractable mathematical analysis.   Therefore the units of these variables are in arbitrary units of $N$-units, $C_a$-units, $D$-units, since they represent various types of cells and thus, they qualitatively rather than quantitatively describe the response of the inflammatory mediators and by-products each represents. Pathogen, $P$, units are more closely related to numbers of pathogens or colony forming units (CFU), but abstract units $P$-units are simply used as well and this population is scaled by $10^6$/cc. More details about the model development can be found in \cite{reynolds2006reduced}.
	
 The $C_a$ variable maintains a small positive background level when the system is in steady state in the presence of no pathogen insult. Following the setup used in \cite{day2010using}, the reference parameter value for $C_a(0)$ is set to $0.125$ and virtual patients have a value that is $\pm25\%$ of the reference value. In addition, five other parameters as well as the initial conditions for $P$ and $C_a$ are made to have differing (positive) values from the reference set to distinguish one virtual patient from another and from the virtual patient having the reference set of parameters.\\

At the end of 168 hours\footnote{168 hours, which is equivalent to one week, is considered to be a clinically relevant time period over which to observe a patient experiencing a severe inflammatory response and when outcome would most likely could be determined. Our results match this time frame, as for example, from Figure 2 one can see that by 168h it is relatively clear to determine to which equilibrium the time courses are converging, whether septic death or the aseptic death region. }, a healthy outcome of a virtual patient's inflammatory response given a particular initial pathogen is defined as the ending state in which  $D=0$, $P=0$, $N=0$, and $C_a$ is at its initial background level.   On the other hand, a patient is considered to have an unhealthy outcome if the ending state has either $P=0$ with all other variable values elevated above their background non-infection steady state values (\textit{aseptic death outcome}) or when $P$ is also elevated along with the other variables (\textit{septic death outcome}). The fact that a virtual patient may, or may not, return to a healthy state depends on the parameter values and initial condition.  When inflammatory levels of a response are deemed excessive (i.e. $N(t)>0.5$ in this work), then intervention becomes necessary to stabilize the patient to his healthy equilibrium (\textit{i.e.} to {\em homeostasis}). In the work here, we consider two virtual patients whose inflammatory response outcomes are septic death and aseptic death, respectively, in the absence of intervention measures. Below, we summarize for each virtual patient the corresponding parameter values used as well as the initial conditions which represents the state of the patient's inflammatory response some  past time when the initial infection occurred, and the inflammatory response has become elevated enough to warrant therapeutic intervention.
\begin{figure}[!ht]
	\begin{center}
		\hspace{-5mm}	\includegraphics[scale=0.5]{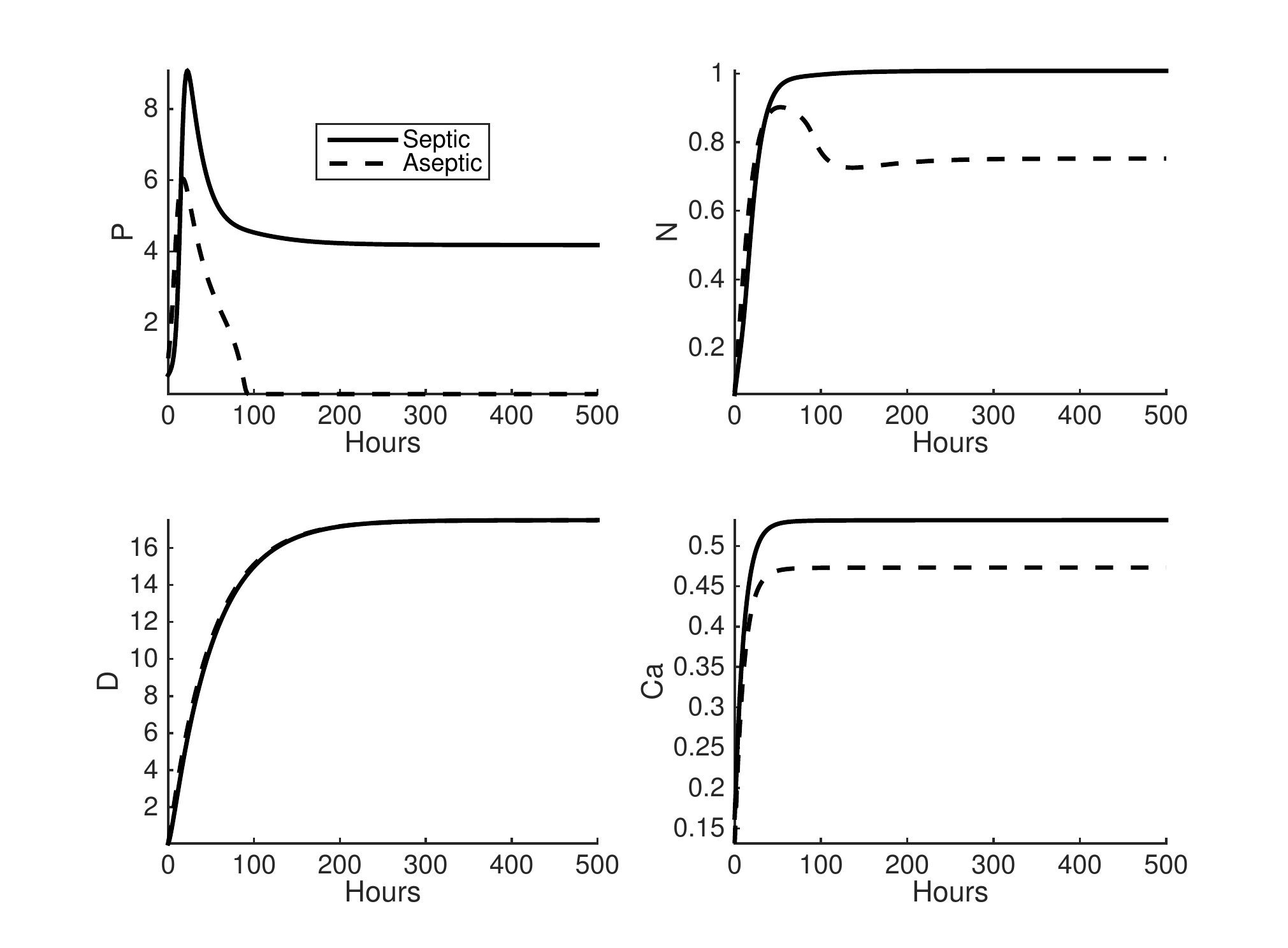}
		\caption{Time courses of system states for natural (open loop) response for virtual patient 1 (solid curves), resulting in a (\emph{septic death outcome}) and for virtual patient 2 (dashed curves), resulting in an (\emph{aseptic death outcome.})}
		\label{Fig:fig2}
	\end{center}
\end{figure}

\begin{itemize}
\item Baseline Patient 
\begin{itemize}
	\item Healthy equilibrium:    $P(0) = 0$, $N(0) = 0$, $D(0) = 0$, $C_a(0) = 0.125$.
	\item Parameter values of the reference patient: $k_{pg} = 0.6$, $k_{cn} = 0.04$, $k_{nd} = 0.02$, $k_{np} = 0.1$, $k_{cnd} = 48$, $k_{nn} =  0.01$.
\end{itemize}
\item Patient 1 (septic death outcome)
\begin{itemize}
	\item Initial conditions:    $P(0) = 0.5360$, $N(0) = 0.0660$, $D(0) = 0.0477$, $C_a(0) = 0.1635$.
	\item Parameter values differing from reference set: $k_{pg} = 0.5846$, $k_{cn} = 0.0409$, $k_{nd} = 0.0242$, $k_{np} = 0.1211$, $k_{cnd} = 49.1243$, $k_{nn} =  0.012$.
\end{itemize}
\item Patient 2 (aseptic death outcome)
\begin{itemize}
	\item Initial conditions:    $P(0) = 1.0017$, $N(0) = 0.0711$, $D(0) = 0.0732$, $C_a(0) = 0.1314$.
	\item Parameter values differing from reference set: $k_{pg} = 0.4746$, $k_{cn} = 0.0386$, $k_{nd} = 0.0223$, $k_{np} = 0.1116$, $k_{cnd} = 46.3367$, $k_{nn} = 0.0112$.
\end{itemize}
\end{itemize}
 The differential equations model using the parameter values and initial conditions corresponding to each virtual patient is numerically solved and the resulting time courses for each are shown in Figure \ref{Fig:fig2}. Solid curves display the solution time courses for patient 1 and dashed curves are for patient 2. As illustrated, in the absence of any control measures, each virtual patient's trajectory evolves toward either of the final unhealthy steady states. Next, we demonstrate how interventions via positive inputs into the system can redirect trajectories away from the unhealthy steady states and toward the healthy steady state.

\section{The Optimal Control Problems} \label{sec:optproblem}
%



We consider optimal control problems in which the  controls  enter linearly into the system. After algebraically rearranging the expressions in equations (\ref{eq:mod1})-(\ref{eq:mod4}) and  renaming the resulting parameter expressions from $a$ to $w$ as shown in Table \ref{tab:tabpar}, the model system with controls can now be written as:
\begin{align} \label{eq:simpmod}
\begin{pmatrix}\dot{P}\\\dot{N}\\\dot{D}\\\dot{C_a}\end{pmatrix}=f(x)+g(x)u= \begin{pmatrix}a P - bP^2 - \frac{cP}{d+eP}- \frac{z}{g+C_a^2}PN  \vspace{2mm} \\
	\frac{hP+ iN+jD}{kP +lN+mD+nC_a^2+o}-pN  \vspace{2mm}\\
	\frac{qN^6}{(r+sC_a^2)^6+tN^6}-\mu_dD \vspace{2mm}\\
	s_c+\frac{uN+vD}{g+C_a^2+gN+wD}-\mu_c C_a  \end{pmatrix}   +\begin{pmatrix}0&0\\1&0\\0&0\\0&1  \end{pmatrix}\begin{pmatrix} u_p\\u_a
\end{pmatrix}
\end{align}

\begin{table}[h!]
\centering
\begin{tabular}{cccccc}
	\toprule
	$a=k_{pg}$  & $e=k_{mp}$ &  $i=c_{\infty}^2s_{nr}k_{nn}$  &  $m=c^2_{\infty}k_{nd}$ &  $q=c^{12}_{\infty}k_{dn}$  \\
	
	$b=\frac{k_{pg}}{P_{\infty}}$ & $z=k_{pn}c_{\infty}^2$ &$j=c^2_{\infty}s_{nr}k_{nd}$ &$s=k_{dn}$& $n=k_{nr}$ \\
	
	$c=k_{pm}s_m$ & $g=c_{\infty}^2$ &  $k=c^2_{\infty}k_{np}$& $o=c^2_{\infty}k_{nr}$&$t=c^{12}_{\infty}$\\
	
	$d=k_m$ & $h=c_{\infty}^2s_{nr}k_{np}$ &  $l=c^2_{\infty}k_{nn}$& $p=k_n$& $u=c^2_{\infty}k_{cn}$\\
	$v=c^2_{\infty}k_{cn}k_{cnd}$ &$r=c^2_{\infty}k_{dn}$ & $w=c^2_{\infty}k_{cnd}$ &\\
	\bottomrule
\end{tabular}
\caption{Parameters of the algebraically rearranged immune model with control terms as given in (\ref{eq:simpmod}).} \label{tab:tabpar}
\end{table}
where $f(x)$ and $g(x)$ represent a vector and a matrix, respectively,  that can be easily inferred from equation (\ref{eq:simpmod}). The state vector is given by $x=(P, N, D, C_a)^T$ and the control vector by $u=(u_p, u_a)^T$. The above alternate representation of the system, after algebraically rearranging the immune model equations, provides a more useful form for the differentiation that will be carried out below in solving the control problem and allows a better understanding of the role of each state in the system.   Two case studies will be  investigated: an optimal control problem with $L^1$-type objective (linear cost term) versus one with an $L^2$-type objective (quadratic cost term) in $J$. \\
For the quadratic case the objective functional to be minimized is therefore given as
\begin{align}
J(u_p,u_a) = \int_0^{t_f}a_1P(t)^2+ a_2N(t)^2+a_3D(t)^2+a_4C_a(t)^2+B_1u_p(t)^2+B_2u_a(t)^2 dt  \label{Eq:quadObj}
\end{align}
For the linear case, the objective functional is given according to 
\begin{align}
J(u_p,u_a) = \int_0^{t_f}a_1P(t)+ a_2N(t)+a_3D(t)+a_4C_a(t)+B_1u_p(t)+B_2u_a(t)dt  \label{Eq:Objlinear}
\end{align}
 The goal here, for both functionals is to determine the optimal control pair $u_p^*$ and $u_a^*$ such that
\begin{align}
J(u_p^*,u_a^*) = \underset{(u_p,u_a)\in \mathbb{U}}\min J(u_p,u_a) \label{eq:cost}
\end{align}

where  the basic control set is 
$$\mathbb{U}=\{ (u_p,u_a) \in  (L^{\infty}(0, t_f) )^2  | \quad 0 \leq u_p(t)  \leq M_p, 0 \leq u_a (t)  \leq M_a, \quad a.e.\}$$
and other state constraints may be considered later.

 Note that during our simulation we also added the final cost term to our objective functional as  $\phi(x(t_f))=P(t_f)+N(t_f)+D(t_f)+C_a(t_f)$ for the linear case and similarly  $\phi(x(t_f))=x(t_f)^Tx(t_f)=P^2(t_f)+N^2(t_f)+D^2(t_f)+C_a^2(t_f)$ for the quadratic one.

The optimal state $x^*$, although satisfying the equality constraints $\dot{x} =f(x(t),u(t))$,  is further constrained by the initial and final boundary condition, or point constraints. In our case, the initial state is $x(0)=(P(0), N(0), D(0), C_a(0))=(0.536,\,0.066,\,0.0477,\,0.1635)$ for patient 1 and $(1.0017, 0.0711, 0.0732, 0.1314)$ for patient 2 and the final state for each is $x(t_f)=(P(t_f), N(t_f), D(t_f), C_a(t_f))=(0,\,0,\,0,\,0.125)$. 


The approach known as the \textit{indirect method}, used in our previous publications \cite{baraoptL1,bara2015}, attempts to preserve the infinite dimensional nature of the problem \cite{subchan2009computational} and uses Pontraygin's maximum principle to solve a TPBVP.  The disadvantage of the indirect method lies in the fact that the boundary problem with the transversality condition is often difficult to solve, and can become further encumbered in the presence of path constraints. In this work, we take an approach that finds a finite dimensional representation necessary to solve the problem. Such approaches are known as \textit{direct methods} and to achieve this, we discretize the state and/or control variables and formulate a nonlinear programming problem (NLP) to solve the optimization problem.

\subsection{Existence of an optimal control}
Practically, before attempting to numerically find an optimal control solution, it is interesting to see if an optimal control  indeed exists. 
The existence of an optimal control for the $L^1$ and $L^2$ objective functional can be obtained from standard results like in \cite{fleming1975,fister1998} or by using a weak convergence direct approach \cite{kelly2016impact}. These results use the bounded range of the controls and the states as well as the structure of the ODE system being linear in the controls. Also the right hand sides of the ODEs are continuous functions of the states.
\begin{asu}\label{asu1}
	\label{asu1i}
	If  the generated inflammatory doses $u_p$ and $u_a$ and the initial states are non-negative  (i.e. $u_p(t) \geq 0, u_a(t) \geq 0  \quad \forall t $ and $P(0)\geq0, \quad N(0)\geq0, \quad D(0)\geq0, \quad C_a(0)\geq0 $) then the solution trajectories of the differential equations (\ref{eq:mod1})-(\ref{eq:mod4}) remain in the positive octant.
\end{asu}

This is a reasonable assumption since the states represent biological entities that are inherently positive.

\section{Numerical solution of optimal control problem}\label{sec:numeric}
\vspace{1mm}
Numerical solutions of our optimal control problem are found by applying direct optimization via the use of PSOPT \cite{becerra2010solving}, an open source pseudospectral optimization software package. PSOPT is written in C++ and uses direct collocation methods, including  Legendre or Chebyshev pseudospectral discretization as well as local transcriptions such as the trapezoidal or Hermite-Simpson integration methods. Sparse nonlinear programming is then used to find local optimal solutions. PSOPT can  be interfaced to, IPOPT, the open source Interior Point Optimization solver (\cite{wachter2006implementation}) for solving large scale optimization problems. This solver will be used for all simulations in this work.

A good example that highlights the relationship between nonlinear programming (NLP) and optimal control can be found in \cite{betts2010practical}, where the Karush-Kuhn-Tucker (KKT) necessary conditions approach the necessary condition of the maximum principle as the number of nodes increases (i.e. as the step size $h\rightarrow0$. The NLP Lagrange multipliers can be interpreted as a discrete approximation of the adjoint variables found using the maximum principle, which motivates the following remark.
\begin{remark}
	Using the Lagrange multipliers provided by PSOPT, the numerical solutions can be verified by showing that they satisfy the  necessary conditions of optimality with good accuracy.
\end{remark}
We observed that the Hermite-Simpson integration method provided better convergence results and accuracy $(\mathcal{O}(h^4))$ (\cite{becerra2010solving}), when compared to either trapezoidal methods provided in PSOPT which are known to have an accuracy of $\mathcal{O}(h^2)$. The Pseudospectral discretization method of Legendre or Chebyshev proved to be slow to converge and not accurate although these have been successfully applied elsewhere (e.g., \cite{bedrossian2009zero}, \cite{ross2004pseudospectral}). In addition, since we found that the solution of the optimal control problem of the previously mentioned two patients to be very similar, we will focus on providing simulation results related to the patient with a septic outcome.

We note that the control objectives in this work are to minimize the levels of pathogen and of damage such that $P=0$ and $D=0$, respectively, while at the same time making sure to use the lowest possible dosing 
amounts. In addition, we make the following assumption
\begin{asu}\label{asu1}
	\label{asu1i}
	If  the generated inflammatory doses $u_p$ and $u_a$ and the initial states are positive semidefinite (i.e. $u_p(t) \geq 0, u_a(t) \geq 0  \quad \forall t $ and $P(0)\geq0, \quad N(0)\geq0, \quad D(0)\geq0, \quad C_a(0)\geq0 $) then the solution trajectories of the differential equations (\ref{eq:mod1})-(\ref{eq:mod4}) remain in the positive octant.
\end{asu}

A number of grid points $N=1000$ is used together with a tolerance of $10^{-6}$ for the NLP solver specifications. All simulations are performed for a week ($168\, hours$) and, therefore, the sampling period, given the number of grid points, is every $0.168$ hours. For each simulation we report the cost function $J_z$, where $z$ specifies the corresponding subsection in which the particular cost function is defined; the entry and exit times, $t_1$ and $t_2$, respectively, of the singular arc, if one exists in the case where an $L^1$ type objective function is used; and the numerical value of the states, $P$, $N$, $D$, and $C_a$ at the final time, $t_f$. Note that in all the control simulations it can be verified that by further integrating the system from the final state $x(t_f)$ the virtual patient will eventually converge to its healthy equilibrium. In the next section we will formulate  the optimal control problem with an $L^1$ type objective and define the necessary optimality condition of the maximum principle related to problem  given in (\ref{Eq:Objlinear}).

\subsection{Optimal Control with $L^1$ type Objective}\label{OC_L1}
We now formulate the optimal control problem with an $L^1$ type objective and define the necessary optimality conditions of Pontraygin's maximum principle. The $L^1$-type objective function has the following form:
\begin{align}
J_{L^1} = {x(t_f)}+\int_0^{t_f}  \Big[{Ax(t)}+ {Bu(t)}\Big] \, dt, \label{Eq:costL1}
\end{align}
where $x=[x_1,x_2,x_3,x_4]^T$ and where $A=[a_1,a_2,a_3,a_4]$ together with $B=[b_1,b_2]$ are row vectors  with positive constants. Also, $x(t)=[P(t),\, N(t),\, D(t),\, C_a]^T \geq 0$ and $u(t)=[u_p, u_a]^T \geq 0$, where $\geq$ means componentwise.\\

 To ensure that the lowest possible dosing amounts are used to achieve the objective, we include the input doses, $u_p$ and $u_a$, in the cost function $J_{L^1}$ above. The mediators  $N$ and $C_a$ are also included in the cost function since the steady state for these values in healthy patients is assumed to be very low for $C_a$ and zero for $N$. This inclusion will ensure to keep the level of inflammatory mediators at lower values. It was also noticed that the optimal control solver converges faster in the presence of $N$ and $C_a$  in the cost function.\footnote{The necessary inclusion here of $N$ and $C_a$ in the objective function differs from that used in our previous studies using MPC (\cite{day2010using}, \cite{zitelli2015combining}).}  \\
 \begin{remark}
	Note that the existence of optimal controls with state constraints are more difficult but can be obtained when variations can be done within the admissible set of controls. For background on problems with state constraints and the corresponding necessary conditions, see \cite{clarke2010optimal,rockafellar1972state} 
 \end{remark}
 The consequence of using an objective that is linear in the controls is that the resulting time-varying controls will be ``bang-bang'', meaning that treatment protocols alternate between the maximum dose of pro- and/or anti-inflammatory therapies and the minimum (i.e. zero) dose. Singular controls may also arise and  will be discussed later. Since we are minimizing the cost of $J_{L^1}$, the standard Hamiltonian is given by:
\begin{align}
H(x,u,\lambda)=Ax(t)+ Bu(t)+\lambda(t)^T(f(x)+g(x)u(t)),\qquad \lambda \in \mathbb{R}^4.
\end{align}

In what follows, we characterize the optimal control which is then illustrated in several simulations where we explore various constraints on the dosing controls as well as mixed state-control constraints.
\begin{theorem} (Characterization of the optimal control)
	Given two optimal controls $u_p^*$, $u_a^*$ with the corresponding solutions $P^*$, $N^*$, $D^*$, $C_a^*$ of the ordinary differential equation given in Eq.(\ref{eq:mod1})-Eq.(\ref{eq:mod4}), there exist adjoint equations associated with the system states given by
 \begin{align}
	\dot \lambda_1 &= -a_1-\lambda_1\Big[a-2bP-\frac{c}{d+eP}+\frac{ceP}{(d+eP)^2}  \nonumber \\[2mm]
	&- \frac{zN}{g+C_a^2}\Big] -\lambda_2\Big[\frac{h(nC_a^2+o)}{(nC_a^2+kP+lN+mD+o)^2}  \Big]; \label{Eq:adj1}  \\[3mm]
	\dot \lambda_2 &= -a_2- \lambda_1\left(\frac{zP}{g+C_a^2}\right)-\lambda_2\Big[\frac{i(nC_a^2+o)}{(nC_a^2+kP+lN+mD+o)^2} \nonumber \\[3mm]
	&-p \Big] -\lambda_3\Big[\frac{6qtN^{11}}{((sC_a^2+r)^6+tN^6)^2}\Big]-\lambda_4\Big[\frac{u(g+C_a^2)}{(g+C_a^2+wD+gN)^2} \Big]; \nonumber \label{Eq:adj2} \\[2mm]
	\dot \lambda_3&= -a_3-\lambda_2\Big[\frac{j(nC_a^2+o)}{(nC_a^2+kP+lN+mD+o)^2} \Big]  \nonumber \\[2mm]
	&+ \lambda_3k_d -\lambda_4\Big[\frac{v(g+C_a^4)}{(g+wD+N)^2} \Big];
	\end{align}
	\begin{align}
	\dot \lambda_4=-a_4-2\lambda_1\frac{zPNC_a}{(g+C_a^2)^2} +\lambda_2\Big[\frac{2nC_a(hP+iN+jD)}{(nC_a^2+kP+lN+mD+o)^2} \Big] \nonumber \\
	+\lambda_3\Big[\frac{12qsN^6(sC_a^2+r)^5}{((sC_a^2+r)^6+tN^6)^2} \Big]
	-\lambda_4\Big[\frac{2(uN+vD)C_a}{(gN+wD+C_a^2+g)^2}-k_c \Big]; \label{Eq:adj4}
	\end{align}
	and satisfying the following transversality condition
	\begin{align}
	\mathbf{\lambda}(t_f)=\Big\{\frac{\partial \phi[{x}(t_f)]}{\partial {x}} \Big\}=1.
	\end{align}
	Furthermore,
	\begin{equation}
	u_p^* = \left\{\begin{array}{lll}
	0 & if & b_1+\lambda_2>0 \\
	N_{max} & if & b_1+\lambda_2<0
	\end{array} \right.
	\end{equation}
	
	\begin{equation}
	u_a^* = \left\{\begin{array}{lll}
	0 & if & b_2+\lambda_4>0 \\
	C_{a_{max}} & if & b_2+\lambda_4<0
	\end{array} \right.
	\end{equation}
\end{theorem} \vspace{3mm}

\begin{proof}
	The adjoint equations and the transversality conditions follow from Pontryagin's maximum principle \cite{pontryagin1987}. Differentiating the negative of the Hamiltonian with respect to the states, $P$, $N$, $D$, and $C_a$, gives the adjoint system in Eq.(\ref{Eq:adj1})-(\ref{Eq:adj4}). According to the maximum principle, the optimality equations are then given by
	\begin{align}
	\frac{\partial H}{\partial u_p} &= 	b_1  +\lambda_2=0 \\
	\frac{\partial H}{\partial u_a} &= b_2  +\lambda_4=0.
	\end{align}
	
	These are known as switching functions and if we denote them by
	\begin{align}
	\phi_1(t)&=	b_1  +\lambda_2 \label{eq:switchfunc1}\textrm{ and}\\
	\phi_2(t)&= b_2  +\lambda_4 ,\label{eq:switchfunc2}
	\end{align}
	then the optimal control that minimizes the Hamiltonian becomes
	\begin{equation}
	u_p^* = \left\{\begin{array}{lll}
	0 & if & \phi_1(t)>0 \\
	N_{max} & if & \phi_1(t)<0\\
	\end{array} \right.
	\end{equation}
	\begin{equation}
	u_a^* = \left\{\begin{array}{lll}
	0 & if & \phi_2(t)>0 \\
	C_{a_{max}} & if & \phi_2(t)<0.\\
	
	\end{array} \right.
	\end{equation}

\end{proof} \vspace{2mm}

  As was done in our previous studies, we define $N_{max}=0.5$ and $C_{a_{max}}=0.62$ as the maximum allowed levels of the variables $N$ and $C_a$, respectively. This restriction is to emulate a realistic bound appropriate in a clinical setting to avoid toxicity effects of potential overdosing.  As will be seen, these maximum bounds are violated when state constraints are not included in the problem statement. We next provide several numerical solutions to the optimal control problem with the $L^1$ type objective and discuss issues regarding when singular controls arise as well as the case necessitating an explicit inclusion of state constraints. The following weights appearing in the objective function \ref{Eq:costL1} for the $L^1$ case will be used for the various constraint scenarios explored:  $(b_1,b_2)=(1,50)$ and $(a_1,a_2,a_3,a_4)=(100, 5, 30, 10)$. The largest weights are assigned to the states $P$ and $D$ since eliminating the pathogen and decreasing the level of damage is our goal.\\

\subsubsection{Numerical results for $L^1$ objective with dosing constraints $ \,0\leq u_p\leq 0.5$ and $0\leq u_a \leq 0.62$}\label{sec:L1dosing}

\begin{figure}[ht]
	\begin{center}
		\includegraphics[scale=0.8]{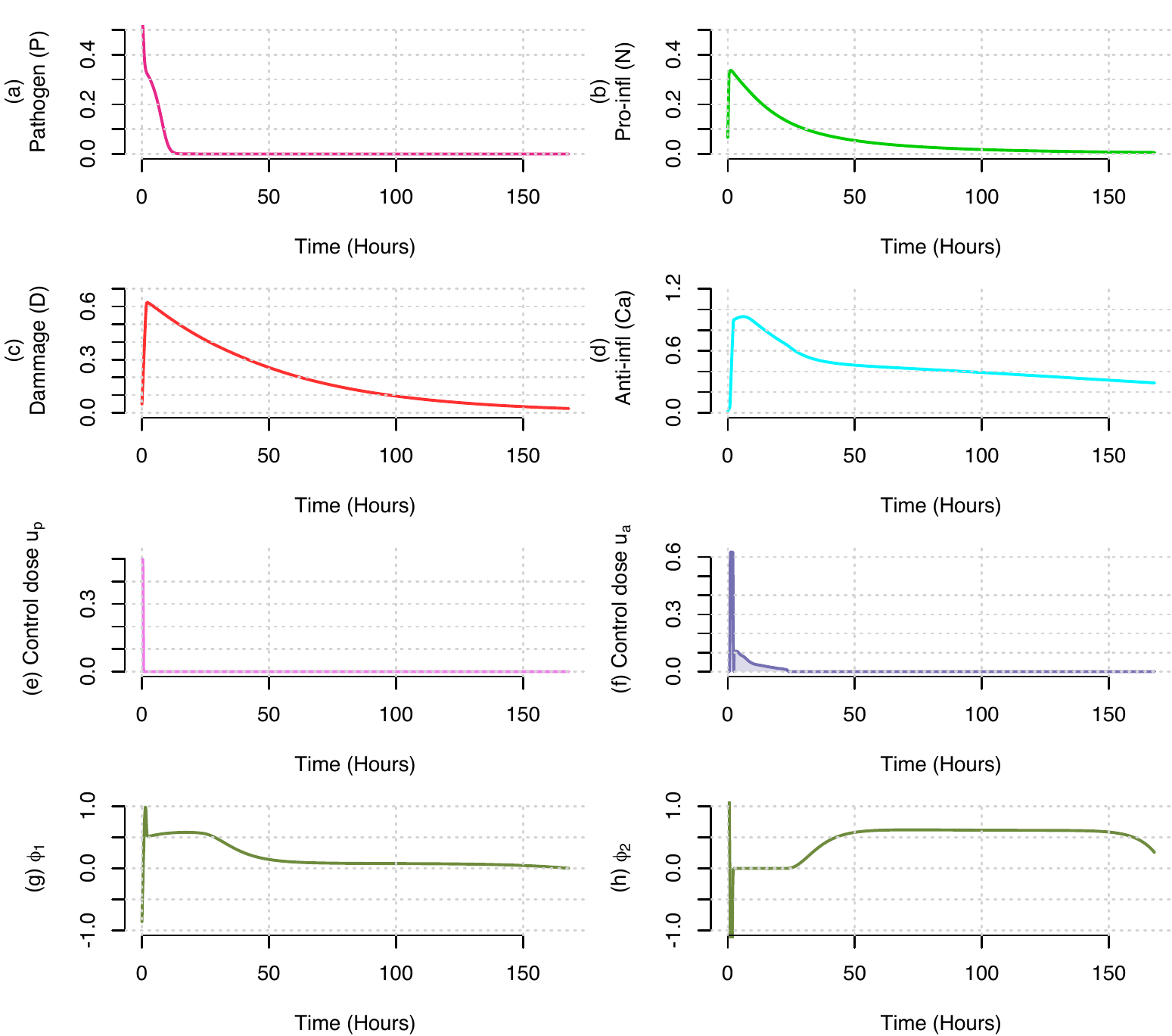}
		\caption{(a)-(d) Time courses for optimal states, (e)-(f) control states, and (g)-(h) the switching functions $\phi_1$ and $\phi_2$. Control constraints in this simulation which uses an $L^1$ type objective are $0\leq u_p\leq 0.5, $ $0\leq u_a \leq 0.62$. }  \label{Fig:bang1}
	\end{center}
\end{figure}

In this scenario we assume the upper bounds on the dosing input to be the same as those used in \cite{day2010using}: $ \,0\leq u_p\leq 0.5$ and $0\leq u_a \leq 0.62$.  Figure \ref{Fig:bang1} displays the simulation results of this scenarios, showing the optimal states of Pathogen ($P$), the pro-inflammatory  state ($N$), the level of Damage ($D$) and the anti-inflammatory state ($C_a$) in Figure \ref{Fig:bang1} (a)-(d) and the resulting control doses, $u_p$ and $u_a$, in panels (e)-(f). The pro-inflammatory control $u_p$ is provided at maximum level for about half an hour, while the anti-inflammatory dose $u_a$ is followed just after for about one hour. Although not very evident in Figure \ref{Fig:bang1}, $u_a$ is applied with one sample difference after $u_p$.\\

When a switching function is not zero there is ``bang-bang'' control and, although the controls are  ``bang-bang'' for most of the simulation, $u_a$ actually portrays a singular control between time $(t_1,t_2)=(2.354,23.88)$. The scaled switching function $\phi_2$ shown in Figure \ref{Fig:bang1}(h) is zero during this interval, unlike $\phi_1$ in panel (g) which does not vanish in any finite open interval. We first demonstrate that the singular control is minimizing and then demonstrate that an approximation to the singular control gives similar results and provides a more practical dosing protocol to implement.

To determine if the singular control generated is minimizing, we derive a solution corresponding to the singular control. However, the fact that there is no explicit dependence on the control in either switching function (\ref{eq:switchfunc1}) or (\ref{eq:switchfunc2}) makes it problematic to determine where they vanish and higher order derivatives must be examined in order to derive a solution that will correspond to the singular control. As mentioned, singular control occurs in the anti-inflammatory control $u_a$ and thus, the switching function corresponding to $\phi_2$ will be of interest. The first derivative of $\phi_2$ is given in Eq.(\ref{eq:phi2}); however, since the control does not appear in the expression, a second derivative is needed.\footnote{If $\phi_2$ vanishes on $(t_1\,,t_2)$ then we must have $\phi_2=\dot{\phi}_2=\ddot{\phi}_2=0$ until the first control appears}

	 \begin{align}
		\dot{\phi_2} =-a_4-2\lambda_1\frac{zPNC_a}{(g+C_a^2)^2}\label{eq:phi2} +\lambda_2\Big[\frac{2nC_a(hP+iN+jD)}{(nC_a^2+kP+lN+mD+o)^2} \Big] \nonumber \\
		+\lambda_3\Big[\frac{12qsN^6(sC_a^2+r)^5}{((sC_a^2+r)^6+tN^6)^2} \Big]
		-\lambda_4\Big[-\frac{2(uN+vD)C_a}{(gN+wD+C_a^2+g)^2}-k_c \Big]=0
		\end{align}

The second derivative of $\phi_2$ is given by Equation (\ref{eq:pphi1}) where the control now appears, and which $\phi_2=0$ can be solved to derive the solution corresponding to the singular control.
\begin{align}
	\ddot\phi_2 &= \left(-\frac{2\lambda_1zNC_a}{(C_a^2+g)^2}-\frac{2b_1nC_ah}{(nC_a^2+kP+lN+mD+o)^2}+\frac{4b_1nC_a(hP+iN+jD)k}{(nC_a^2+kP+lN+mD+o)^3}\right)\dot{P}\nonumber\\
	 &+\theta_1\left(\frac{(hP+iN+jD)}{(nC_a^2+kP+lN+mD+o)}-pN\right)+\left[\frac{-2b_1nC_aj}{(nC_a^2+kP+lN+mD+o)^2}\nonumber \right.\\
	 &\left. + \frac{4b_1nC_a(hP+iN+jD)m}{(nC_a^2+kP+lN+mD+o)^3}+b_2\left(\frac{2vC_a}{(gN+wD+C_a^2+g)^2}-\frac{2(2uN+2vD)C_aw}{(gN+wD+C_a^2+g)^3}\right)\right]\dot{D}\nonumber\\
	 &+\theta_2\left(\frac{sc+(uN+vD)}{(gN+wD+C_a^2+g)}-k_cC_a\right)-2\dot{\lambda}_1\frac{zPNC_a}{(g+C_a^2)^2}+\dot{\lambda}_2\Big[\frac{2nC_a(hP+iN+jD)}{(nC_a^2+kP+lN+mD+o)^2} \Big]\nonumber \\
	&+ \theta_1u_p+\theta_2u_a=0, \label{eq:pphi1}
\end{align}
where
\begin{align}
\theta_1 &= -2\frac{\lambda_1zPC_a}{(C_a^2+g)^2}-\frac{2b_1niC_a}{(nC_a^2+kP+lN+mD+o)^2}+\frac{4b_1nC_a(hP+iN+jD)l}{(nC_a^2+kP+lN+mD+o)^3}\nonumber\\
&+\frac{72qsN^5(sC_a^2+r)^5}{((sC_a^2+r)^6+tN^6)^2}\lambda_3-\frac{144\lambda_3qsN^{11}(sC_a^2+r)^5t}{((sC_a^2+r)^6+tN^6)^3}\nonumber \\
&+b_2\left(\frac{2uC_a}{(gN+wD+C_a^2+g)^2}-\frac{2(2uN+2vD))C_ag}{(gN+wD+C_a^2+g)^3}\right) \label{Eq:theta1}
\end{align}

and
\begin{align}
\theta_2	 &=8\frac{\lambda_1zPNC_a^2}{(C_a^2+g)^3}-2\frac{\lambda_1fPN}{(C_a^2+g)^2}-8\frac{\lambda_2(hP+iN+jD)n^2C_a^2}{(nC_a^2+kP+lN+mD+o)^3}\nonumber \\
	&+2\frac{\lambda_2(hP+iN+jD)n}{(nC_a^2+kP+lN+mD+o)^2} +120\lambda_3\frac{qs^2N^6(sC_a^2+r)^4C_a}{((sC_a^2+r)^6+tN^6)^2} \nonumber \\
	&-\lambda_4\left(\frac{-(2uN+2vD)}{(gN+wD+C_a^2+g)^2} + \frac{(4(2uN+2vD))C_a^2}{(gN+w*D+C_a^2+g)^3}\right)-288\lambda_3\frac{qs^2N^6(sC_a^2+r)^{10}C_a}{((sC_a^2+r)^6+tN^6)^3}. \label{Eq:theta2}
\end{align}

For singular control to be minimizing, the generalized Legendre-Clebsch condition shown in (\ref{eq:LC}) needs to be satisfied \cite{krener1977high}.
	 \begin{align}
	 	(-1)\frac{\partial}{\partial u_a}\frac{d^{2}}{dt^{2}}\frac{\partial H}{\partial u_a}\geq 0 \label{eq:LC}
	 \end{align}

That is, during the simulation, the expression in $-\theta_2$ (See Figure \ref{Fig:bang1}(h)) must be positive semi-definite. Although not reported here, a plot of $-\theta_2$ (See Eq.(\ref{Eq:theta2})) shows that it is strictly positive for $t \in (t_1,t_2)$, where $t_1=2.354$ and $t_2=23.88$, proving the optimality of the singular arc. If we consider the pro-inflammatory input, $u_p$, to be constant then the singular control part of $u_a$ can be easily derived from Eq.(\ref{eq:pphi1}). Note that in all simulations, $u_p$ is  zero during the interval where a  singular control is present in $u_a$.\\

Figure \ref{Fig:bang1}(a)-(f) displays results of the simulation. The cost function $J_{\ref{sec:L1dosing}}$ associated with Figure \ref{Fig:bang1}, the entry and exit times of the singular control, and the final states corresponding to the simulation are given as follows:

\begin{align*}
	J_{\ref{sec:L1dosing}} &= 2104.78
	  &   t_1& =2.354  &  t_2&=23.88 	\\
	P(t_f)&=0 &    N(t_f)&=0.0061  &  D(t_f)&=0.025 \\
	C_a(t_f)&=0.29
\end{align*}

Finally, we consider an approximation to the ``bang-bang-singular'' control (or what we term the \textit{basic optimal control}), where we replace the time course for the singular control with a constant input, namely zero. In other words, we do not provide any dosing of $u_a$ during the period where the singular control occurred and determine if this approximation is acceptable in terms of the resulting outcomes of the states. The motivation of replacing the singular control with a constant is to suggest a more practical implementation of the dosing protocol.

\begin{figure}[ht]
	\includegraphics[scale=0.66]{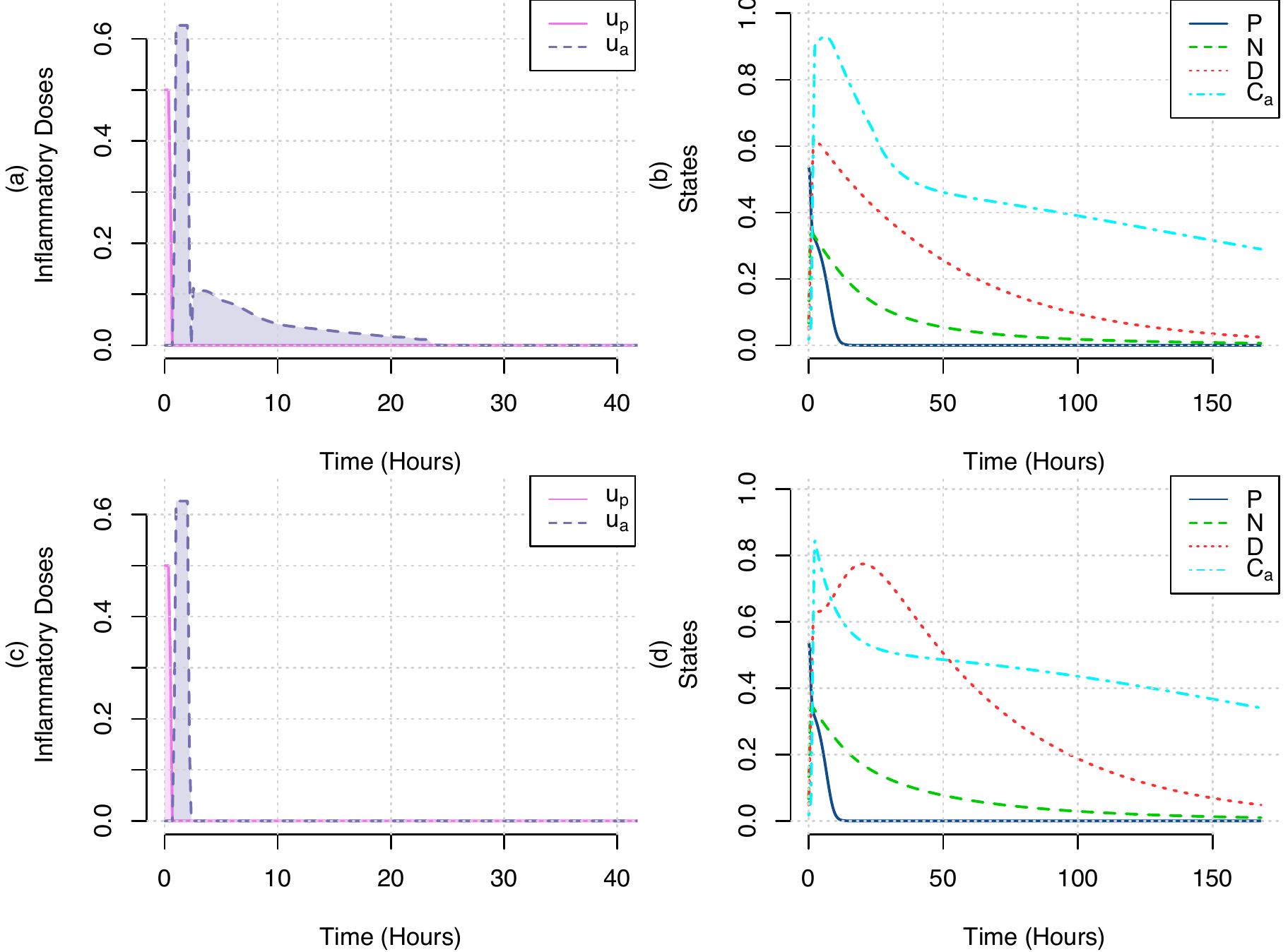}
	\caption{Comparison of basic optimal control with an approximate control for the $L^1$ objective case with dosing constraints of $ \,0\leq u_p\leq 0.5$ and $0\leq u_a \leq 0.62$. (a) Generated optimal controls $u_p$ and  $u_a$  and (b) corresponding optimal states. (c) Generated optimal control $u_p$ and approximation of optimal control $u_a$ and (d) corresponding states. The time frame displayed in both plots (a) and (c) of the control doses is for $t_f=40$ hours since afterward the values of both controls are zero.}  \label{Fig:bangsingular}
\end{figure}

 Figure \ref{Fig:bangsingular} provides the comparison of the simulation using the optimal controls generated for $u_p$ and $u_a$ (panels (a)-(b) which are the same as in Figure \ref{Fig:bang1}(a)-(f)) versus the simulation where an approximation of the $u_a$ control was made for the singular arc (panels (c)-(d)). Comparison of Figure \ref{Fig:bangsingular} panel (d) with (b) shows that the states $C_a$ and $D$ are the ones mainly affected by this approximate control case. In other words, changing the level of anti-inflammatory dose $C_a$, namely reducing it to zero for that time period, has a direct impact on the mitigation of damage ($D$). Thus, a lack of anti-inflammatory mediator at critical moments results in an increase of the level of damage,  which is seen in Figure \ref{Fig:bangsingular}(d). In this particular case, the long-term dynamics of the two simulations are the same; however, such approximations would need to be done carefully. An improved solution might be to approximate the singular arc by a piece-wise constant function of say, four pieces which would better capture the gradual decrease of $u_a$ occurring in the singular arc.

\subsubsection{Numerical results for $L^1$ objective with dosing constraints $0\leq u_p\leq 0.25$ and $0\leq u_a \leq 0.31$}\label{sec:L1halvedcon}
\begin{figure}[h!]
 	\includegraphics[scale=0.65]{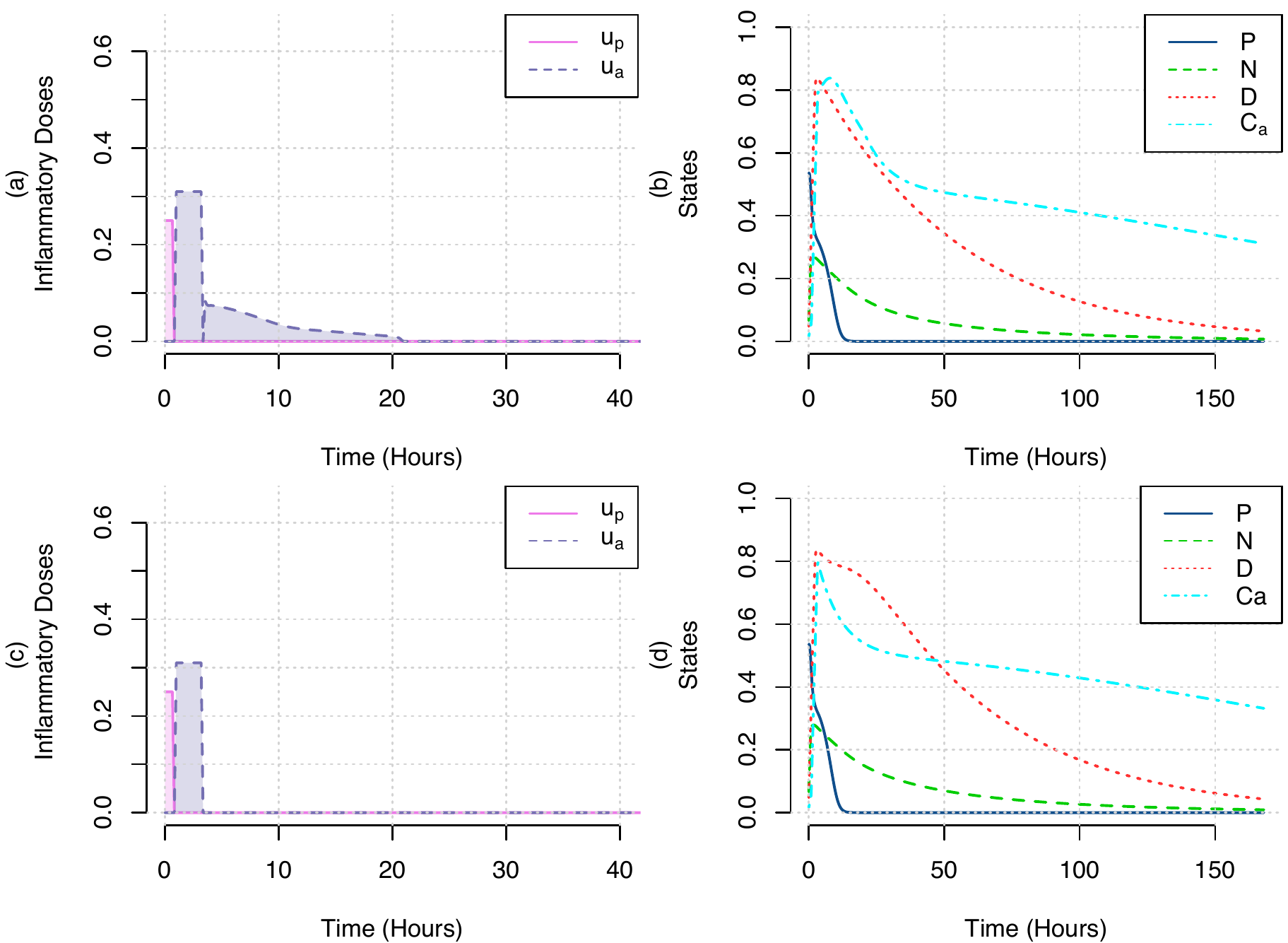}
 	\caption{Comparison of basic optimal control with an approximate control for the $L^1$ objective case with upper bounds on dosing constraints given as $u_{pmax}=0.25$ and $u_{amax}=0.31$. (a) Generated optimal controls $u_p$ and  $u_a$  and (b) corresponding optimal states. (c) Generated optimal control $u_p$ and approximation of optimal control $u_a$ and (d) corresponding states. The time frame displayed in both plots (a) and (c) of the control doses is for $t_f=40$ hours since afterward the values of both controls are zero.}  \label{Fig:bang2}
 \end{figure}

We next explored how a reduction of the maximum allowable doses affected the ability of the control problem with $L^1$ objective to drive the patient to its healthy state. Thus the upper bounds on $u_p$ and $u_a$ were reduced to $0\leq u_p\leq 0.25$ and $0\leq u_a \leq 0.31$, which is half of the amount used in the previous section. The resulting simulation is given in Figure \ref{Fig:bang2} (a)-(b). 
The numerical information about the cost, exit and entry times for the singular arc, as well as the final states corresponding to Figure \ref{Fig:bang2} are given by

\begin{align*}
J_{\ref{sec:L1halvedcon}} &= 2469.904
&   t_1& =3.363  &  t_2&=21.02 	\\
P(t_f)&=0 &    N(t_f)&=0.0074  &  D(t_f)&=0.0326 \\
C_a(t_f)&=0.311
\end{align*}
 The singular control verifies the Legendre-Clebsch condition in (\ref{eq:LC}) for optimality. The simulation results show an  increase in the cost function $J_{\ref{sec:L1halvedcon}} $ when compared to $J_{\ref{sec:L1dosing}} $ of the previous section and this can be explained by  the increase of the level of Damage, $D$, as a result of decreasing the total dose of the anti-inflammatory mediator, $C_a$. Similar to what was done in the previous section, Figure \ref{Fig:bang2} (d) displays the outcome of the states when an approximation of the optimal control input, $u_a$, is made as shown in panel (c).\\

  In either sets of simulations for this and previous section and in either the optimal or approximate cases, we observe that the level of damage ($D$) remains such that $D\leq 1$ (and can be verified to approach zero eventually in the absence of further control) and the level of pathogen ($P$) decreases relatively quickly to zero. An interesting pattern we observe is that each dose associated with an increase of the pro-inflammatory mediator is followed immediately after by a dose of anti-inflammatory therapy. This behavior can be explained by the fact that the immune system requires an initial boost of activated phagocytes, $N$, to eliminate the pathogen threat, but the resulting inflammation causes some self-harm seen in the relative increase of tissue damage, $D$, which then needs the inhibitory effects of $C_A$ via an anti-inflammatory dose, $u_a$, to ensure the decay of damage.\\ 
  
  The anti-inflammatory dosing input is applied for a period of time clearly longer than $u_p$. It is interesting to notice in Figure \ref{Fig:bang2} that when the maximum allowable doses were decreased by half, the control doses seems to be  applied for a longer duration when compared to the initial scenario where $N_{max}=0.5$ and $C_{amax}=0.62$. As a consequence of lowering the intensity of the dose, the optimal control compensates by allowing a longer application. This behavior was also observed when solving the TPBVP for the  $L^1$ problem using the forward backward sweep method in \cite{lenhart2007optimal}.\\

  It can also be seen in Figure \ref{Fig:bangsingular} of the previous section that the anti-inflammatory states do not remain below the upper bound mentioned previously on the anti-inflammatory state $C_a$ (i.e. $C_{amax}=0.62$\footnote{The max value for $C_a$ was derived from the bifurcation study in \cite{reynolds2006reduced} of the $N/D$ subsystem. It was determined in that study that the system loses bi-stability between the healthy and the death state if $C_a>0.62$ with only the healthy state being stable. For this reason, to ensure the presence of a stable death state the latter bound on $C_a$ is enforced. The max value of $N=0.5$ was determined experimentally (through simulations; See \cite{day2010using})}) and the same is used for the simulations of this section given in Figure \ref{Fig:bang2}, even though the maximum upper bound on the doses for both controls was reduced by half. Thus, in the next section, we impose state constraints on $N$ and $C_a$ in the problem formulation.\\

\subsubsection{Numerical results for $L^1$ objective with state and dosing constraints}\label{sec:L2bothcon}
In this section, we now impose the state constraints on $N$ and $C_a$ we defined earlier of $0\leq N(t)\leq 0.5$ and $0\leq C_a\leq 0.62$, along with dosing constraints of $0\leq u_p\leq 0.5$ and $0\leq u_a \leq 0.62$ which are the same dosing constraints of Section \ref{sec:L1halvedcon}. The analysis of the numerical results obtained when using the state constraints together with the control constraints would be subject, in most part, to the same discussion as in the section before; the only difference would be the additional jump conditions on the adjoints and the transversality condition. To be more pertinent, we will refer the reader to the similar discussion presented later in section \ref{sec:L2pure}. For the sake of brevity, the details are not provided for this case. 

\begin{figure}[ht]
	\includegraphics[scale=0.7]{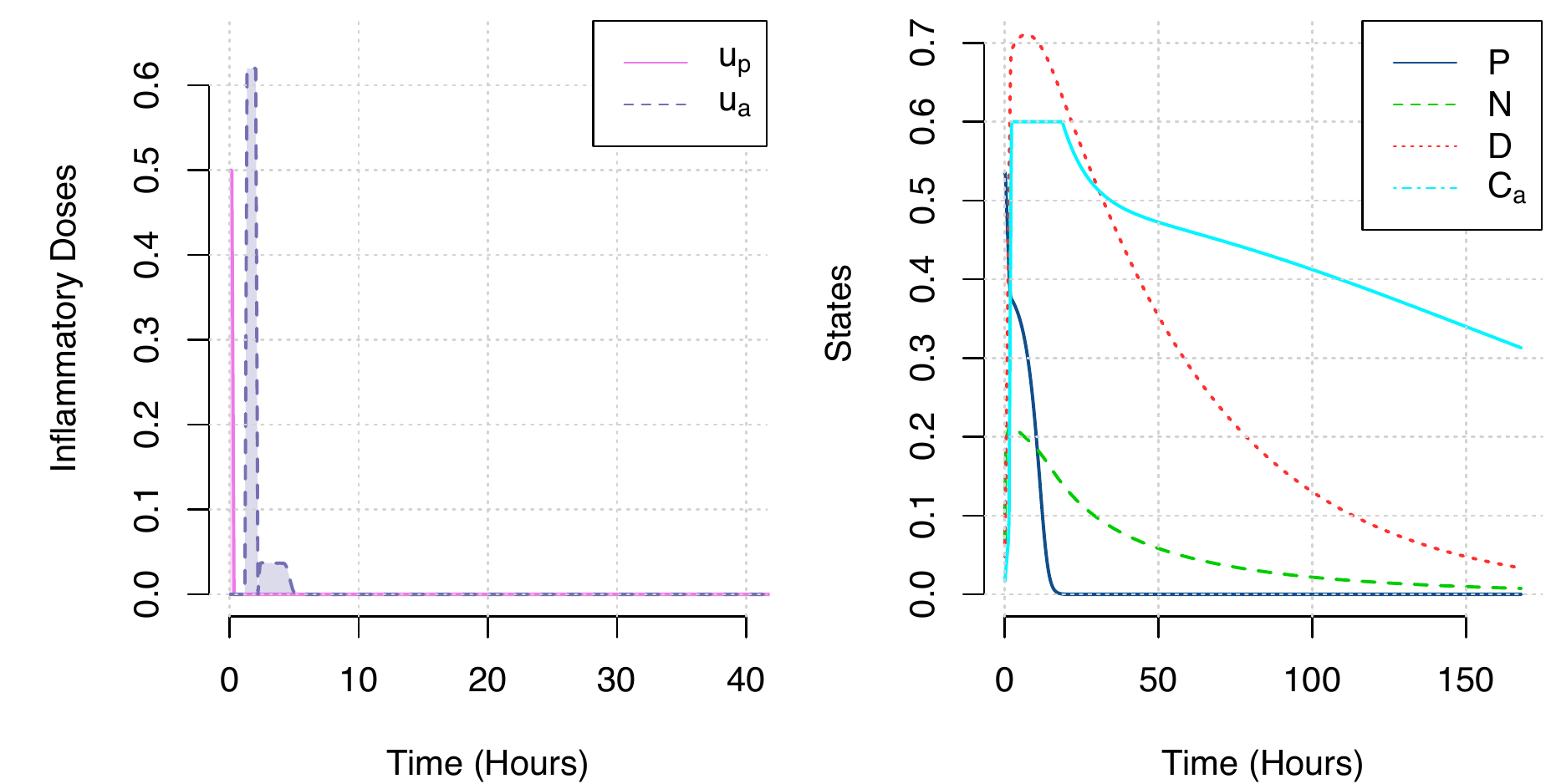}
	\caption{ Optimal doses and corresponding state trajectories for an $L^1$ type objective with state constraints $0\leq N(t)\leq 0.5$, $0\leq C_a\leq 0.62$ together with dosing constraints $0\leq u_p\leq 0.5$ and $0\leq u_a \leq 0.62$. (a) Dosing profiles for $u_p$ and $u_a$ controls that lead to (b) a healthy outcome as all states decrease toward background levels. The time frame displayed in (a) of the control doses is for $t_f=40$ hours since afterward the values of both controls are zero.}  \label{Fig:bang31}
\end{figure}

The resulting simulation is given in Figure \ref{Fig:bang31} and the numerical information about the costs and final states corresponding to Figure \ref{Fig:bang31} are given by

\begin{align*}
J_{\ref{sec:L2bothcon}} &= 2393.339
&   t_1& =2.186  &  t_2&=4.036 	\\
P(t_f)&=0.0000 &    N(t_f)&=0.0073  &  D(t_f)&=0.0318 \\
C_a(t_f)&=0.3093
\end{align*}

The simulation in Figure \ref{Fig:bang31}(a) shows that a singular arc occurs between time $t_1=2.186$ and $t_2=4.036$, which can be verified as before by checking the resulting function $\phi_2$, as the one given in Figure \ref{Fig:bang1}(h). The optimality of the singular arc can be derived by verifying the Legendre-Clebsh condition of (\ref{eq:LC}). Notice that \ref{Fig:bang31}(a) shows that the pro-inflammatory dose is applied for a duration of almost 20 minutes followed by a rest period of 50 minutes before the anti-inflammatory dose is injected for a longer duration corresponding to one hour.\\

Thus far, we have studied optimal control scenarios in which the objective is linear in the control. In the remainder of this paper the focus will be on deriving the optimal control doses subject to an objective that is quadratic in the control as well as different constraint types. A discussion comparing the results of the uses of the different objective functions will also be given.

\subsection{Optimal Control with $L^2$ type Objective}
In this section we present numerical simulations of the optimal control problem when the objective is quadratic in the control.
\begin{align}
J_{L^2} = x^T(t_f)x(t_f)  +   \int_0^{t_f} \Big[x(t)^TAx(t) + u(t)^TBu(t) \Big] dt \label{Eq:costL2}
\end{align}
where A and B are matrices with diagonals $(a_1,a_2,a_3,a_4)$ and $B=(b_1,b_2)$ respectively. Also, $x(t)=[P(t),\, N(t),\, D(t),\, C_a]^T \geq 0$ and $u(t)=[u_p, u_a]^T \geq 0$, where $\geq$ means componentwise.\\

We also once again define $N_{max}=0.5$ and $C_{a_{max}}=0.62$ as the maximum allowed levels of the variables $N$ and $C_a$, respectively. We provide several numerical solutions to the optimal control problem with the $L^2$ type objective under various constraint definitions. The following weights appearing in the objective function (\ref{Eq:costL2}) for the $L^2$ case will be used for the various constraint scenarios explored in the subsections that follow: $(b_1,b_2)=(1,20)$ and $(a_1,a_2,a_3,a_4)=(10, 1, 10, 10)$. Note that the objective function includes the cost at the final time,  $\phi(x(t_f)^2)$.  All other aspects of the numerical settings are as given in the beginning of Section \ref{sec:numeric}.\\

The Hamiltonian for the $L^2$ type objective is given according to
\begin{align}
H(x,u,\lambda)={x(t)}^TAx(t)+ {u(t)}^TBu(t)+\lambda(t)^T(f(x)+g(x)u) \label{eq:hamiltQ}
\end{align}
where $\lambda(t)=(\lambda_1, \lambda_2, \lambda_3, \lambda_4)^T$. While The $L^1$ type objective leads in general to more intricate analysis involving bang-bang and singular control, the $L^2$ type objective analysis often results in continuous control, when no state constraints are included. The case of mixed state-control inequality constraint will be introduced in next section, followed by a section that looks at separate or \emph{pure} control and state constraints. For all constraints used in the following subsections, those associated with $N$ and/or $u_p$ have a maximum bound of $0.5$, while those associated with $C_a$ and/or $u_a$ have a maximum bound of $0.62$.

\subsubsection{Numerical results for $L^2$ objective with mixed control-state inequality constraints}\label{sec:L2mixed}

\begin{figure}[ht]
	\hspace{-4mm}\includegraphics[scale=0.7]{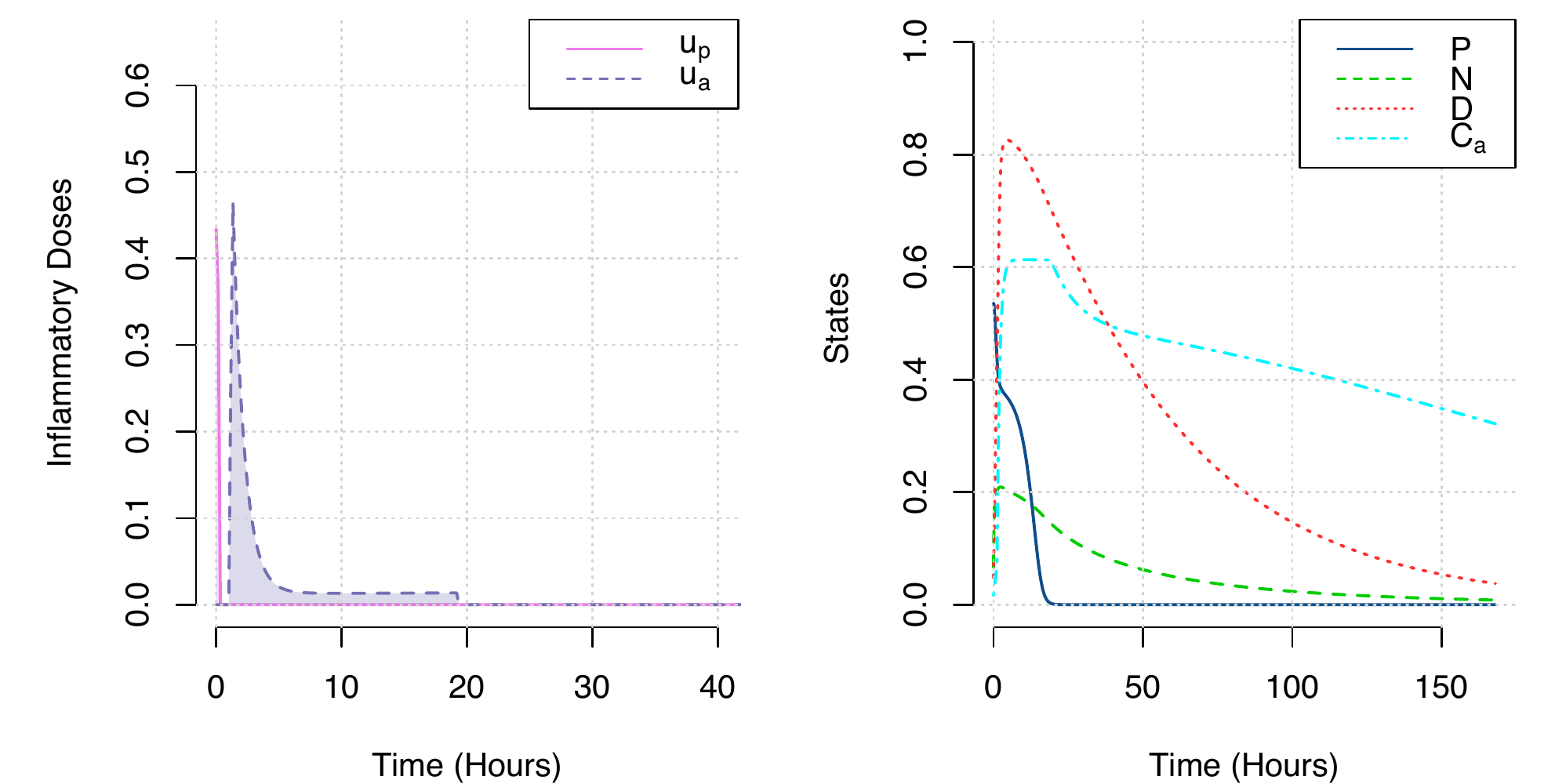}
	\caption{(a) Optimal doses and (b) corresponding state trajectories for a quadratic type objective with mixed control state inequality constraints $0\leq u_p(t) \leq 0.5-N(t)$, $0\leq u_a(t) \leq 0.62-C_a(t)$. The time frame displayed in (a) of the control doses is for $t_f=40$ hours since afterward the values of both controls are zero.}  \label{Fig:mixedcon}
\end{figure}

Consider the mixed control-state inequality constraints \footnote{Note that we adhere to the sign convention in \cite{bryson} } given according to
\begin{align}
C_1(N(t),u_p(t)) &=u_p(t)+N(t)-0.5 \leq  0 \label{eq:C1}\\
C_2(C_a(t),u_a(t)) &=u_a(t)+C_a(t)-0.62 \leq 0 \label{eq:C2}.
\end{align}
Imposing mixed control-sate constraints allows the derived optimal control to account for the actual value of the states, and hence only the necessary control dose will be administered.  
The simulation results for this OCP are given in Figure \ref{Fig:mixedcon} and correspond to the numerical information about the costs and final states given as follows:
\begin{align}\label{eq:costmixed}
J_{\ref{sec:L2mixed}} &= 602.67&   P(t_f)&=0 & N(t_f)&=0.0082 & D(t_f)&=0.0376 & C_a(t_f)&=0.3216
\end{align}

Note that both controls are constrained to be positive semidefinite, i.e., $0\leq u_p$ and $0\leq u_a$. All the constraints are adjoined such that the augmented Hamiltonian in equation (\ref{eq:hamiltQ}) can be expressed according to

\begin{align}\label{eq:augham}
\mathcal{H}(x,u,\lambda,\mu) &= H(x,u,\lambda)+\mu_1 C_1(N,u_p)+\mu_2C_2(C_a,u_a)- \mu_3u_p-\mu_4u_a.
\end{align}
Necessary conditions for minimizing the Hamiltonian can then be derived. The multiplier $\mu$ is given such that it verifies the complementary condition
\begin{align}
\mu_i = \begin{cases} 0 & if\quad C_i<0\\ \mu_i \geq 0 & if\quad C_i=0\quad for\quad i=1,2.	\end{cases}
\end{align}
The adjoint  equations are given in the following compact form
\begin{align}
\lambda^T=-\mathcal{H}_x=\begin{cases} -H_x(x,u,\lambda) & if \quad C<0\\-H_x(x,u,\lambda)-\mu C_x & if \quad C=0;	\end{cases} \label{adj:mixed}
\end{align}
and more explicitly as follows:
\begin{align}
\dot\lambda_1 &= -\frac{\partial H}{\partial P} & \dot\lambda_3 &= -\frac{\partial H}{\partial D}\\
\dot\lambda_2 &= -\frac{\partial H}{\partial N}-\mu_1 & \dot\lambda_4 &= -\frac{\partial H}{\partial C_a}-\mu_2.
\end{align}
The controls $u_p$ or $u_a$ can be derived from the mixed constraints along the constrained arc. In other words, whenever
\begin{align}
C(x,u)&=0 & for\quad t\in(t_1,t_2).
\end{align}
Hence, the boundary controls are given by
\begin{align}
u_{p_b} &= 0.5-N(t) \label{Eq:mix1}\\
u_{a_b} & =0.62-C_a(t). \label{Eq:mix2}
\end{align}
The free controls can be derived from the minimum condition such that
\begin{align}
u_p &= -\frac{\lambda_2}{2b_1}\\
u_a & = -\frac{\lambda_4}{2b_2}.
\end{align}
Since the regularity condition is verified, i.e., $C_u(x,u)\neq0$, the multiplier $\mu$ follows from the local minimum condition $\mu=\mu(x,\lambda)=-H_u(x,\lambda,u_b(x))/C_u(x,u_b(x))$. In other words, by
\begin{align}
\mu_1&= -2b_1u_p-\lambda_2 \label{Eq:mu1}\\
\mu_2&=-2b_2u_a-\lambda_4, \label{Eq:mu2}
\end{align}
whenever the lower constraints on the controls are not active (i.e. $u_p>0$ and $u_a>0$). When the mixed constraints are active, the multipliers $\mu_1$ and $\mu_2$ are given by substituting the equation for $u_{p_b}$ given in (\ref{Eq:mix1}) into (\ref{Eq:mu1}) and similarly for $u_{a_b}$ of (\ref{Eq:mix2}) into (\ref{Eq:mu2}), such that
 \begin{align}
 \mu_1&= -2b_1(0.5-N(t))-\lambda_2 \\
 \mu_2&=-2b_2(0.62-C_a(t))-\lambda_4.
 \end{align}
The optimal controls are given as follows
\begin{align}
u_p(t) = \begin{cases} -\frac{1}{2b_1}(\lambda_2+\mu_1) & if\quad 0\leq t \leq t_0=1.009\\ 0 & if\quad t_0\leq t\leq t_f=168;	\end{cases}
\end{align}
\begin{align}
u_a(t) = \begin{cases} 0, & if\quad 0\leq t\leq t_1 = 1.009\\
\frac{-1}{2b_2}(\lambda_4+\mu_2), & if \quad t_1\leq t \leq t_2=19\\
 0, & if\quad t_2\leq t\leq t_f=168.	\end{cases}
\end{align}
Note that the mixed control-state constraints, $C_1$ and $C_2$ in equations (\ref{eq:C1}) and (\ref{eq:C2}), respectively, are active whenever the controls are not zero; that is, $\mu_1\geq 0$ and $\mu_2\geq 0$. For $0\leq t\leq t_0$ and $t_1\leq t\leq t_2$ the switching structure is given by the junction conditions
\begin{align}
\mu_1(t)&= -2b_1(0.5-N(t))-\lambda_2(t) =0\\
\mu_2(t)&=-2b_2(0.62-C_a(t))-\lambda_4(t)=0.
\end{align}

The next section presents a scenario where we have pure state constraints on $N$ and $C_a$ as well as constraints on the upper bounds of the control doses.

\subsubsection{Numerical results for $L^2$ objective with pure state and control constraints}\label{sec:L2pure}

\begin{figure}[!ht]
	\hspace{-2mm}\includegraphics[scale=0.7]{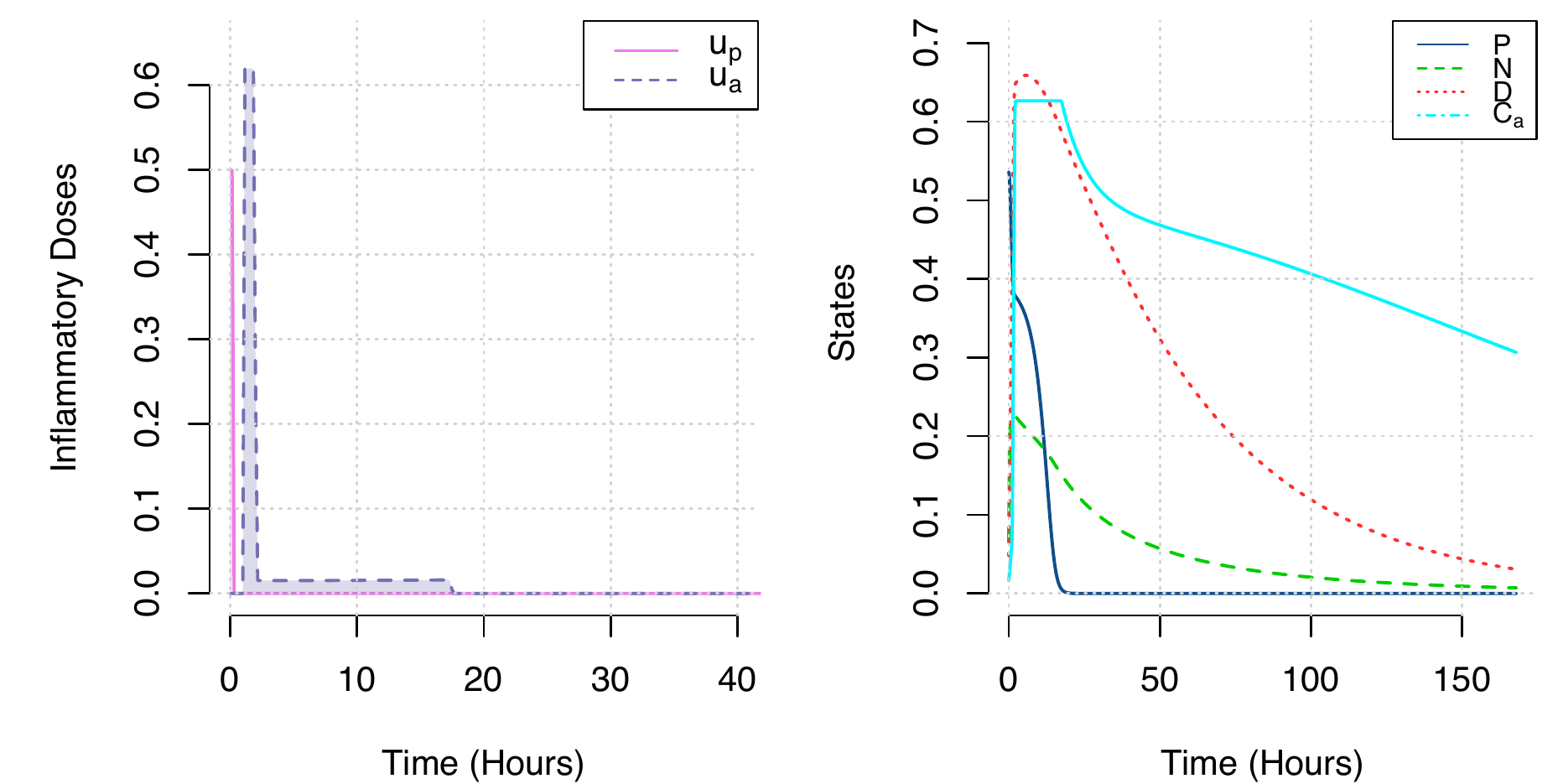}
	\caption{ (a) Optimal doses and (b) corresponding state trajectories for a quadratic type objective with state path constraints $0\leq N(t)\leq 0.5$, $0\leq C_a\leq 0.62$. The time frame displayed in (a) of the control doses is for $t_f=40$ hours since afterward the values of both controls are zero.}  \label{Fig:mixed}
\end{figure}

Now consider imposing path constraints on the states $N$ and $C_a$ such that $N(t)\leq 0.5$ and $C_a(t)\leq 0.62$ $\forall t$. In addition, we will require as before that the following control constraints be met: $0\leq u_p\leq 0.5$ and $0\leq u_a\leq 0.62$. The simulation results for this OCP are given in Figure \ref{Fig:mixed} and correspond to the numerical information about the costs and final states given as follows:

\begin{align*}
J_{\ref{sec:L2pure}} &= 513.77 &  P(t_f)&=0 &    N(t_f)&=0.0071   &D(t_f)&=0.0306 & C_a(t_f)&=0.3065
\end{align*}

Consider the last inequality in the (OCP) problem $S(x(t)) \leq 0$, given specifically for this case as
\begin{align}
S_1(N(t)) &= N(t)-0.5 \leq  0 \\
S_2(C_a(t)) &= C_a(t)-0.62 \leq 0.
\end{align}

Whenever one of the constraints becomes active on a subinterval we have
\begin{align}
S_i(x(t))&=0 &\text{for all} &&t\in(t_1,t_2)\subset [0,168].
\end{align}

Let $r$ be the \emph{order} of the state constraint; that is, the smallest non-negative numbers necessary for the control to first appear \cite{sethi2000optimal}. It is straightforward to see that $r=1$ and furthermore,
\begin{align}
\dot{S}_2(x,u_a) = \dot{C}_a = s_c+\frac{uN+vD}{g+C_a^2+gN+wD}-\mu_c C_a +u_a(t). \label{eq:S1}
\end{align}
The state constraint satisfies the following regularity condition
\begin{align}
\frac{\partial}{\partial u_a}\dot{S}_2(x,u_a)=(0,1)\neq(0,0).
\end{align}

\begin{remark}
	We can derive similar conditions for $S_1$. We will focus on the state $C_a$ since that is the only state that becomes active.
\end{remark}

The boundary control $u_{a_b}$ on the constrained arc can be derived from (\ref{eq:S1}) by also replacing $C_a(t)=0.62$ such that
\begin{align}
u_{a_b}=-s_c-\frac{uN+vD}{g+0.62^2+gN+wD}+0.62\mu_c.
\end{align}
 The augmented Hamiltonian is given according to

\begin{align}
\mathcal{H}(x,u,\lambda,\lambda_s)&={x(t)}^TAx(t)+ {u(t)}^TBu(t)+\lambda(t)^T(f(x)+g(x)u)+\lambda_s S(x),\\
&+ \mu_1 (u_p(t)-0.5)+\mu_2(u_a(t)-0.62)-\mu_3u_p-\mu_4u_a,
\end{align}
where $\lambda_s = [\lambda_{sp}, \,\lambda_{sa}]$ is a row vector and $S(x)=[S_1,\,S_2]^T$ is a column vector.
The adjoint equation can be formulated in a similar fashion as in (\ref{adj:mixed}). The multiplier $\lambda_{sa}$ needs to verify the complementary slackness condition \cite{sethi2000optimal} given as:
\begin{align}
\lambda_{sa} = \begin{cases}
0      & if \quad  S_2<0 \\
\lambda_{sa}\geq 0 & if \quad S_2 =0. \label{Eq:complementary}
\end{cases}
\end{align}

\begin{figure}[ht]
	\hspace{-2mm}\includegraphics[scale=0.65]{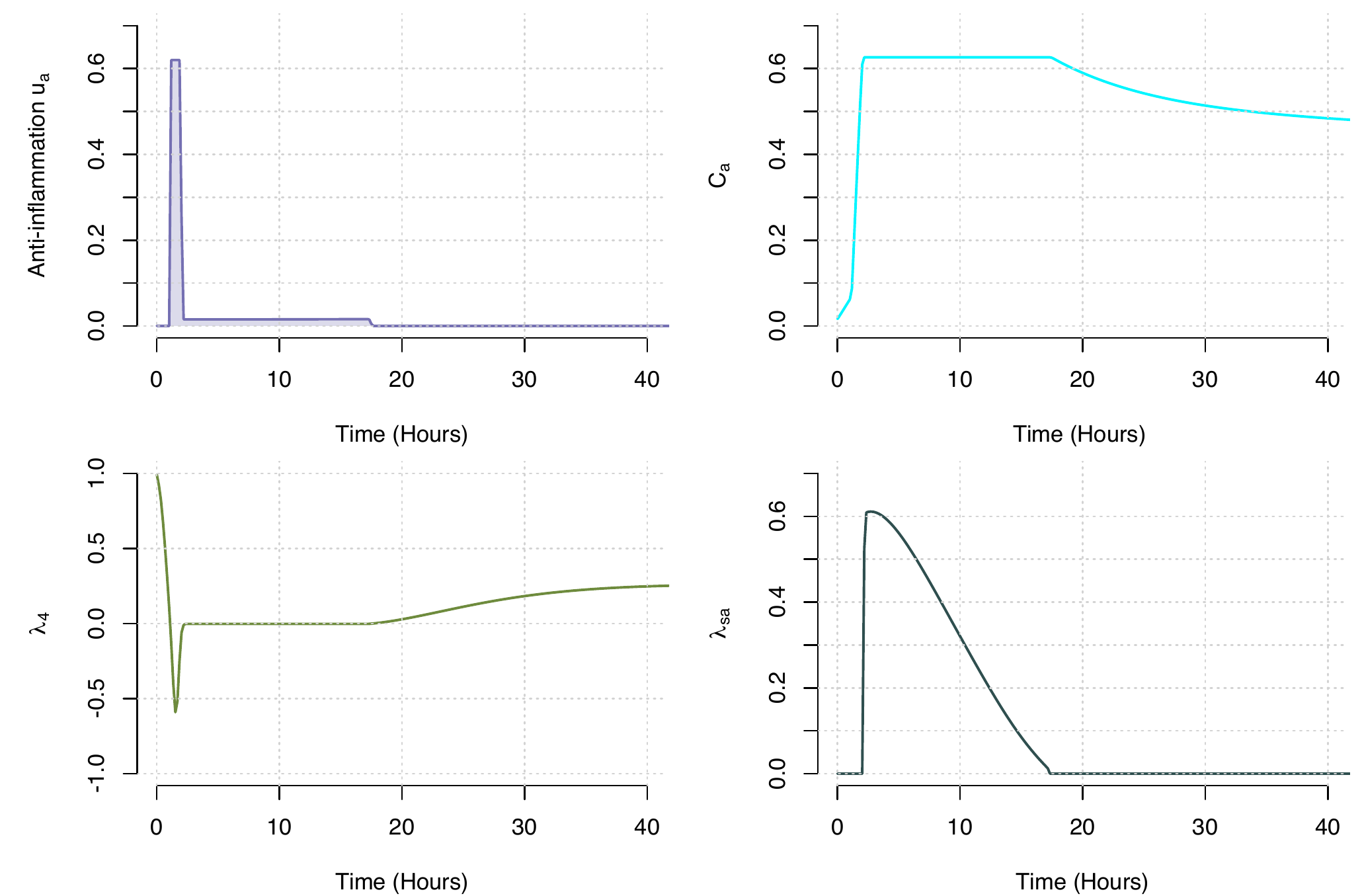}
	\caption{Numerical solutions of the state and control constrained problem with $L^2$ objective. Top row: control $u_a(t)$ and constrained anti-inflammatory state $C_a(t)$. Bottom row: The scaled adjoint $\lambda_4$ corresponding to the state $C_a(t)$ together with the multiplier $\lambda_{sa}$ corresponding to the constraint $S(x(t))=C_a(t)-0.62$. The time frame displayed in the figures is for $t_f=40$ hours since afterward the values of the control inputs are zero. }  \label{Fig:mulambda}
\end{figure}

The adjoint equation corresponding to the dynamic equation $\dot{C_a}$ is given by
\begin{align}
\dot\lambda_4(t) = -\frac{\partial \mathcal{H}}{\partial C_a}-\lambda_{sa}. \label{Eq:end}
\end{align}
The minimum  condition can be derived as in the previous section when $0=\frac{\partial\mathcal{H}}{\partial u_p}=2b_2u_a+\lambda_4$. On the boundary arc, $u_a$ is constant (See Figure \ref{Fig:mulambda}), and we get $\dot{\lambda}_4=0$ as can be seen in the same figure. The multipliers $\lambda_{sa}$ can be derived from  (\ref{Eq:end}) according to
\begin{align}
	\lambda_{sa}=-\frac{\partial \mathcal{H}}{\partial C_a}\Bigg|_{C_a=0.62}.
\end{align}

The optimal controls are given as follows 
\begin{align}
u_p(t) = \begin{cases} 0.5 & if\quad 0\leq t \leq t_0=0.3363\\
 0 & if\quad t_0\leq t\leq t_f=168;	\end{cases}
\end{align}
\begin{align}
u_a(t) = \begin{cases} 0, & if\quad 0\leq t\leq t_1 = 1.009\\
0.62, & if \quad t_1\leq t \leq t_2=2.186 \\
u_{a_b}, & if \quad t_2\leq t \leq t_3=17.32 \\
\frac{-1}{2b_2}(\lambda_4), & if \quad t_3\leq t \leq t_4=17.66\\
0, & if\quad t_4\leq t\leq t_f=168.	\end{cases}
\end{align}
The jump conditions and transversality condition are given, respectively, as
\begin{align}
\lambda(t_{s^{-}})&=\lambda(t_{s^{+}})+\beta_sS_x(x(t_s))\quad \beta_s\geq 0 \textrm{ and }\label{Eq:jump}\\
\lambda(t_f)&=(2P(t_f), 2N(t_f), 2D(t_f), 2C_a(t_f)).
\end{align}

Figure \ref{Fig:mulambda} displays the computed anti-inflammatory dose $u_a$ and the evolution of the constrained state $C_a$  together with the adjoint variable $\lambda_4$ and the multiplier $\lambda_s$ corresponding to the state constraint $S(x(t))=C_a(t)-0.62$. The Figure outlines that the constraint on $C_a$ becomes active in the interval $[2.186,17.32]$. One can also see in the same figure that the multiplier $\lambda_s$ satisfies the complementary condition (\ref{Eq:complementary}). The jump discontinuity at $t_1=2.186$ on the control $u_a$ for this $L^2$ optimal control is nonintuitive since one would expect a continuous control as in the scenario of the previous section; however, this result is common when using pure state constraints (See Eq. (\ref{Eq:jump})).
The scaled adjoint variable $\lambda_4$ in Figure \ref{Fig:mulambda}, which is associated with the differential equation $\dot{x}_4$, has a jump discontinuity at the time the constraint becomes active.
\begin{remark}
The objective functional value corresponding to section \ref{sec:L2mixed}, which is associated with the mixed control-state constraints is higher than the one given in section \ref{sec:L2pure}, that is associated with pure control and state constraints. In the case that imposes mixed control-state inequality constraints, the maximum allowed doses are updated based on the the actual states of $C_a$ and $N$. For instance, if at a sampling time, $t$, $C_a(t)=0.62$ the upper bound on $u_a$ is zero and similarly, if $N(t)=0.5$, then the upper bound on $u_p$ is zero. However, this constraint formulation is more restrictive with the amount of doses given which results in an increase in the damage response (Figure \ref{Fig:mixedcon}) compared with the pure state and control constraints (Figure \ref{Fig:mixed}), which may require additional or longer dosing to achieve a successful outcome.
\end{remark}

\section{Conclusion} \label{sec:conclusion}

In this work, we determined successful treatment protocols to regulate the inflammatory response modeled with an ordinary differential equations system through solving the optimal control problem using a direct approach.  Two cases were considered regarding the form of the objective functional:  linear versus quadratic in the control. The results of the direct numerical methods used are in agreement with the characterization of the optimal control. A noticeable difference was not observed in the way the controls are provided when considering an $L^1$ versus $L^2$ type objective.  In either case, it was observed that state constraints needed to be included explicitly in order for the states not to exceed the prescribed bounds set for realistic considerations. The case where singular control appears is also addressed  when the objective is linear in the control, the optimality of the singular arcs was demonstrated by verifying the Legendre-Clebsch condition. 

The dynamics of inflammatory immune response is more complex and varied than represented by the reduced inflammatory model considered in this study, however the obtained computational results can provide insight into possible treatment strategies and the methods could prove to be useful tools to incorporate in future practice. In particular, one can see a prevailing pattern in the way the doses are implemented in this study and is in agreement with results of previous computational studies (\cite{day2010using}, \cite{zitelli2015combining}).  Whether considering the solution for the $L^2$ or $L^1$ type objective, generally speaking, a successful control strategy is marked by an initial bolus of a pro-inflammatory control ($u_p$) administered first to quickly eliminate the harmful effect of the pathogen, followed just after by a two phase application of the anti-inflammatory mediator ($u_a$). In the first phase, a bolus of anti-inflammation is applied, followed by a longer duration application of $u_a$, but with a very small amount. The similar patterns in the generated control doses $u_p$ and $u_a$ when solving for either type of objective seems clinically reasonable when thinking that in the presence of a pathogen, the pro-inflammatory mediators will naturally increase to try to eliminate the pathogen, which in turn will cause a certain increase in the level of collateral tissue damage which incites further inflammation. If not controlled, this positive feedback may result in the death of the patient. Fortunately, however, the response also includes a release of anti-inflammatory mediators which helps alleviate the harmful effect of the pro-inflammation, ultimately driving the patient to its healthy equilibrium in the absence of further pathogen threat. So essentially, the dynamic necessary for stabilizing a patient with sever inflammation is reflected exactly by increasing the pro-inflammatory effects with input $u_p$ followed just after by a dose of anti-inflammatory mediator, with input $u_a$. It is important to notice in our study that we assume the dosages represent the actual concentrations with idealized effect on patient; that is, we do not assume any model for the pharmacokinetics or pharmacodynamics. The effects of these two phenomena on optimal protocols could be included in further study.  Future work would also benefit from the use of an immune response model that displays similar complex behavior as the reduced model but does not lump too many variables together to see how the obtained results would then relate to what has been presented in this work.


\begin{thebibliography}{10}
	\expandafter\ifx\csname url\endcsname\relax
	\def\url#1{\texttt{#1}}\fi
	\expandafter\ifx\csname urlprefix\endcsname\relax\def\urlprefix{URL }\fi
	\expandafter\ifx\csname href\endcsname\relax
	\def\href#1#2{#2} \def\path#1{#1}\fi
	
	\bibitem{cohen2002immunopathogenesis}
	J.~Cohen, The immunopathogenesis of sepsis, Nature 420~(6917) (2002) 885--891.
	
	\bibitem{FDAXigris}
	U.~Food, D.~A. (FDA)'',
	\href{http://www.fda.gov/Drugs/DrugSafety/ucm277114.htm}{Fda drug safety
		communication: Voluntary market withdrawal of xigris [drotrecogin alfa
		(activated)] due to failure to show a survival benefit} (Oct. 2011).
	\newline\urlprefix\url{http://www.fda.gov/Drugs/DrugSafety/ucm277114.htm}
	
	\bibitem{Grau1997}
	G.~E. Grau, D.~N. Maennel, Tnf inhibition and sepsis � sounding a cautionary
	note, Nature Medicine 3 (1997) 1193--1195.
	
	\bibitem{Reinhart2001}
	K.~Reinhart, W.~Karzai, Anti-tumor necrosis factor therapy in sepsis: update on
	clinical trials and lessons learned, Crit Care Med 29~(7 Supplement 1) (2001)
	S121--S125.
	
	\bibitem{parker2010systems}
	R.~S. Parker, G.~Clermont, Systems engineering medicine: engineering the
	inflammation response to infectious and traumatic challenges, Journal of The
	Royal Society Interface (2010) rsif20090517.
	
	\bibitem{bequette2012challenges}
	B.~W. Bequette, Challenges and recent progress in the development of a
	closed-loop artificial pancreas, Annual reviews in control 36~(2) (2012)
	255--266.
	
	\bibitem{doyle2014closed}
	F.~J. Doyle, L.~M. Huyett, J.~B. Lee, H.~C. Zisser, E.~Dassau, Closed-loop
	artificial pancreas systems: engineering the algorithms, Diabetes care 37~(5)
	(2014) 1191--1197.
	
	\bibitem{haddad2003}
	W.~M. Haddad, T.~Hayakawa, J.~M. Bailey, Adaptive control for non-negative and
	compartmental dynamical systems with applications to general anesthesia, Int.
	J. Adapt. Control Signal Process 17 (2003) 209--235.
	
	\bibitem{stengel2004stochastic}
	R.~F. Stengel, R.~Ghigliazza, Stochastic optimal therapy for enhanced immune
	response, Mathematical biosciences 191~(2) (2004) 123--142.
	
	\bibitem{kirschner1997}
	D.~Kirschner, S.~Lenhart, S.~Serbin, Optimal control of the chemotherapy of
	{HIV}, Journal of mathematical biology 35~(7) (1997) 775--792.
	
	\bibitem{stengel2002optimal}
	R.~F. Stengel, R.~Ghigliazza, N.~Kulkarni, O.~Laplace, Optimal control of
	innate immune response, Optimal Control Applications and Methods 23~(2)
	(2002) 91--104.
	
	\bibitem{day2010using}
	J.~Day, J.~Rubin, G.~Clermont, Using nonlinear model predictive control to find
	optimal therapeutic strategies to modulate inflammation, Mathematical
	Biosciences and Engineering (MBE) 7~(4) (2010) 739--763.
	
	\bibitem{zitelli2015combining}
	G.~Zitelli, S.~M. Djouadi, J.~D. Day, Combining robust state estimation with
	nonlinear model predictive control to regulate the acute inflammatory
	response to pathogen., Mathematical biosciences and engineering: MBE 12~(5)
	(2015) 1127--1139.
	
	\bibitem{parker2000advanced}
	R.~Parker, E.~Gatzke, F.~Doye, Advanced model predictive control (mpc) for type
	i diabetic patient blood glucose control, in: American Control Conference,
	2000. Proceedings of the 2000, Vol.~5, IEEE, 2000, pp. 3483--3487.
	
	\bibitem{grosman2010zone}
	B.~Grosman, E.~Dassau, H.~C. Zisser, L.~Jovanovi{\v{c}}, F.~J. Doyle, Zone
	model predictive control: a strategy to minimize hyper-and hypoglycemic
	events, Journal of diabetes science and technology 4~(4) (2010) 961--975.
	
	\bibitem{li2013}
	H.~Li, W.~M. Haddad, Model predictive control for a multicompartment
	respiratory system, IEEE Transactions on Control Systems Technology 21~(5)
	(2013) 1988--1995.
	
	\bibitem{bara2015}
	O.~Bara, J.~Day, S.~M. Djouadi, Optimal control of an inflammatory immune
	response model, in: "Proceedings of the IEEE 54th Annual Conference on
	Decision and Control (CDC), IEEE, 2015, pp. 1283--1288.
	\newblock \href {http://dx.doi.org/10.1109/CDC.2015.7402388}
	{\path{doi:10.1109/CDC.2015.7402388}}.
	
	\bibitem{baraoptL1}
	O.~Bara, S.~Djouadi, J.~Day, Immune therapy using optimal control with $l^1$-
	type objective, in: Proceedings of the Americain Control Conference (ACC),
	IEEE, 2016, pp. 4895--4900.
	\newblock \href {http://dx.doi.org/10.1109/ACC.2016.7526128}
	{\path{doi:10.1109/ACC.2016.7526128}}.
	
	\bibitem{reynolds2006reduced}
	A.~Reynolds, J.~Rubin, G.~Clermont, J.~Day, Y.~Vodovotz, G.~B. Ermentrout, A
	reduced mathematical model of the acute inflammatory response: I. derivation
	of model and analysis of anti-inflammation, Journal of theoretical biology
	242~(1) (2006) 220--236.
	
	\bibitem{day2006reduced}
	J.~Day, J.~Rubin, Y.~Vodovotz, C.~C. Chow, A.~Reynolds, G.~Clermont, A reduced
	mathematical model of the acute inflammatory response ii. capturing scenarios
	of repeated endotoxin administration, Journal of theoretical biology 242~(1)
	(2006) 237--256.
	
	\bibitem{pontryagin1987}
	L.~Pontryagin, Mathematical Theory of Optimal Processes, Classics of Soviet
	Mathematics, Taylor \& Francis, 1987.
	
	\bibitem{lenhart2007optimal}
	S.~Lenhart, J.~T. Workman, Optimal control applied to biological models, CRC
	Press, 2007.
	
	\bibitem{bryson}
	A.~E. Bryson, Applied optimal control: optimization, estimation and control,
	CRC Press, 1975.
	
	\bibitem{ledzewicz2007optimal}
	U.~Ledzewicz, H.~Sch{\"a}ttler, Optimal controls for a model with
	pharmacokinetics maximizing bone marrow in cancer chemotherapy, Mathematical
	biosciences 206~(2) (2007) 320--342.
	
	\bibitem{swierniak2003optimal}
	A.~Swierniak, U.~Ledzewicz, H.~Schattler, Optimal control for a class of
	compartmental models in cancer chemotherapy, International Journal of Applied
	Mathematics and Computer Science 13~(3) (2003) 357--368.
	
	\bibitem{subchan2009computational}
	S.~Subchan, R.~Zbikowski, Computational optimal control: tools and practice,
	John Wiley \& Sons, 2009.
	
	\bibitem{fleming1975}
	W.~Fleming, R.~Rishel, Deterministic and stochastic optimal control, Springer,
	1975.
	
	\bibitem{fister1998}
	K.~R. Fister, S.~Lenhart, J.~S. McNally, Optimizing chemotherapy in an hiv
	model, Electronic Journal of Differential Equations 1998~(32) (1998) 1--12.
	
	\bibitem{kelly2016impact}
	M.~R. Kelly~Jr, J.~H. Tien, M.~C. Eisenberg, S.~Lenhart, The impact of spatial
	arrangements on epidemic disease dynamics and intervention strategies,
	Journal of biological dynamics 10~(1) (2016) 222--249.
	
	\bibitem{becerra2010solving}
	V.~M. Becerra, Solving complex optimal control problems at no cost with psopt,
	in: Computer-Aided Control System Design (CACSD), 2010 IEEE International
	Symposium on, IEEE, 2010, pp. 1391--1396.
	
	\bibitem{wachter2006implementation}
	A.~W{\"a}chter, L.~T. Biegler, On the implementation of an interior-point
	filter line-search algorithm for large-scale nonlinear programming,
	Mathematical programming 106~(1) (2006) 25--57.
	
	\bibitem{betts2010practical}
	J.~T. Betts, Practical methods for optimal control and estimation using
	nonlinear programming, Vol.~19, Siam, 2010.
	
	\bibitem{bedrossian2009zero}
	N.~S. Bedrossian, S.~Bhatt, W.~Kang, I.~M. Ross, Zero-propellant maneuver
	guidance, Control Systems, IEEE 29~(5) (2009) 53--73.
	
	\bibitem{ross2004pseudospectral}
	I.~M. Ross, F.~Fahroo, Pseudospectral knotting methods for solving nonsmooth
	optimal control problems, Journal of Guidance, Control, and Dynamics 27~(3)
	(2004) 397--405.
	
	\bibitem{clarke2010optimal}
	F.~Clarke, M.~De~Pinho, Optimal control problems with mixed constraints, SIAM
	Journal on Control and Optimization 48~(7) (2010) 4500--4524.
	
	\bibitem{rockafellar1972state}
	R.~T. Rockafellar, State constraints in convex control problems of bolza, SIAM
	journal on Control 10~(4) (1972) 691--715.
	
	\bibitem{krener1977high}
	A.~J. Krener, The high order maximal principle and its application to singular
	extremals, SIAM Journal on Control and Optimization 15~(2) (1977) 256--293.
	
	\bibitem{sethi2000optimal}
	S.~P. Sethi, G.~L. Thompson, Optimal control theory: applications to management
	science and economics, Springer, 2005.
	
\end{thebibliography}

\end{document}